\documentclass[11pt]{amsart}
\usepackage[a4paper,margin=28mm]{geometry}
\usepackage{amsmath,amssymb,mathrsfs,booktabs,microtype,graphicx}
\usepackage{tikz}
\usepackage[hidelinks]{hyperref}

\numberwithin{equation}{section}
\newtheorem{theorem}{Theorem}[section]
\newtheorem{proposition}[theorem]{Proposition}
\newtheorem{lemma}[theorem]{Lemma}
\newtheorem{corollary}[theorem]{Corollary}
\newtheorem{maintheorem}{Theorem}

\newtheorem{maincorollary}[maintheorem]{Corollary}
\theoremstyle{definition}
\newtheorem{definition}[theorem]{Definition}
\newtheorem{remark}[theorem]{Remark}
\newtheorem{example}[theorem]{Example}

\newcommand{\Mbar}{\overline{\mathcal M}}
\newcommand{\ev}{\mathrm{ev}}
\newcommand{\He}{\mathrm{He}}
\newcommand{\taut}{\mathrm{taut}}
\DeclareMathOperator{\Var}{Var}
\DeclareMathOperator{\Aut}{Aut}
\DeclareMathOperator{\Tr}{Tr}
\DeclareMathOperator{\rk}{rank}

\title{Asymptotic distributions of the Betti and Hodge numbers of $\overline{\mathcal M}_{g,n}$}
\author{Jinwon Choi}
\address{Department of Mathematics and Research Institute of Natural Sciences, Sookmyung Women's University, Seoul 04310, Korea}
\email{jwchoi@sookmyung.ac.kr}

\author{Young-Hoon Kiem}
\address{School of Mathematics, Korea Institute for Advanced Study, 85 Hoegiro, Dongdaemun-gu, Seoul 02455, Korea}
\email{kiem@kias.re.kr}
\date{October 2026}

\begin{document}
\begin{abstract}
We study the asymptotic distributions of the Betti and Hodge numbers of the moduli spaces $\overline{\mathcal M}_{g,n}$ of stable pointed curves and the all-light Hassett spaces $\overline{\mathcal M}_{g,0|n}$, for fixed genus $g\ge1$ as the number of markings tends to infinity. We prove central limit theorems for the full Betti distributions, whose variances grow linearly in $n$ with leading coefficients independent of the genus. For the normalized even Betti numbers, we establish refined local limit theorems with estimates uniform in degree. These local estimates yield asymptotic formulas for individual even Betti numbers and determine the sharp central ultra-log-concavity thresholds $\frac{3(e-2)}{2(3-e)}$ and $3$, respectively. For $g\ge2$, the off-diagonal Hodge components form proportions of exact orders $n^{-1/2}$ and $n^{-1}$ of the total cohomology of the two spaces, respectively. The conditional distributions of the Hodge difference $p-q$ given $p\ne q$ converge, for both spaces, to a single discrete law independent of the genus. We also prove that tautological classes account for all but a proportion of exact order $n^{-1}$ of the even cohomology of $\overline{\mathcal M}_{g,n}$. The tautological degree distribution converges to the standard Gaussian under the centering and scaling of the ambient even Betti distribution. An explicit asymptotic formula for the ratio of the ordinary Euler characteristics of the light and heavy spaces answers a question of Kannan, Serpente and Yun. The proofs combine the heavy/light substitution with a stratification of the light spaces by stable cores, singularity analysis of Hodge generating functions, and refined Fourier estimates exploiting Poincar\'e duality.
\end{abstract}
\maketitle
\section{Introduction}\label{sec:introduction}

The moduli spaces $\Mbar_{g,n}$ of stable $n$-pointed curves of genus $g$ \cite{DM,Knu} are central objects of algebraic geometry, and their cohomology is the natural setting for intersection theory on moduli spaces, Gromov--Witten theory and quantum cohomology \cite{Manin}. Much is known about this cohomology in special ranges. In genus zero, Keel described the cohomology ring and gave recursions for the Betti numbers \cite{Keel}. In genus one, the even cohomology is tautological by a theorem of Petersen \cite{Petersen}, and the Hodge--Deligne polynomials can be computed explicitly \cite{Kiem}. In arbitrary genus, the cohomology in low degrees has been analyzed in recent work of Bergstr\"om, Faber and Payne \cite{BFP} and of Canning, Larson and Payne \cite{CLP}. We refer to \cite{FaberPandharipande} for a survey of tautological and non-tautological classes.

In this paper we take a complementary, global point of view: we fix the genus and let the number of markings grow. The total Betti number of $\Mbar_{g,n}$ grows like $n!\,(e-2)^{-n}$ times a power of $n$. Exact computation of individual Betti numbers is difficult, so we instead seek asymptotic estimates and ask about the \emph{shape} of the cohomology:
\begin{enumerate}
\item how the Betti numbers are distributed over the degrees, and how large an individual Betti number is;
\item what proportion of the cohomology is tautological, and what proportion lies off the Hodge diagonal;
\item which log-concavity properties the Betti numbers satisfy near the middle degree.
\end{enumerate}
Questions of this kind are classical for combinatorial sequences, where Gaussian limit laws are ubiquitous \cite{FS}. In earlier work we established Gaussian limit laws and central ultra-log-concavity for $\Mbar_{0,n}$ and $\mathbb P^1[n]$ \cite{CKPNAS}, and for the Fulton--MacPherson compactifications of configuration spaces of smooth projective curves \cite{CKTJM}. Here we treat all three questions in positive genus.

Throughout, $g\ge1$ is fixed, cohomology has rational coefficients, and implicit constants may depend on $g$. Along with the space $\Mbar_{g,n}$, whose markings have weight one and are distinct, we study the all-light Hassett space
\[
 \Mbar_{g,0|n}=\Mbar_{g,(1/n,\ldots,1/n)}
\]
of \cite{Hassett}, whose markings have weight $1/n$ and may therefore coincide at smooth points of the curve. We call $\Mbar_{g,n}$ the \emph{heavy} space and $\Mbar_{g,0|n}$ the \emph{light} space. Both are smooth and proper Deligne--Mumford stacks of dimension
\[
 d_n=n+3g-3,
\]
so their Betti numbers are palindromic by Poincar\'e duality. We study the light spaces using their geometric stratification and apply the heavy/light substitution of Kannan, Serpente and Yun \cite{KSY} to obtain results for $\Mbar_{g,n}$.

\subsection{Gaussian behavior of the Betti numbers}\label{subsec:intro-gaussian}
Let $X^H_{g,n}$ denote the random variable that takes the value $k$ with probability proportional to $b_{2k}(\Mbar_{g,n})$, and define $X^L_{g,n}$ in the same way using $\Mbar_{g,0|n}$. Thus $X^H_{g,n}$ is half the degree of an even-degree class chosen uniformly at random from a homogeneous basis. By Poincar\'e duality, both variables have mean exactly $d_n/2$. Let
\[
 \nu=\frac{3-e}{6(e-2)}=0.0653685\ldots,
\]
and let $$\Phi(x)=\int_{-\infty}^x\frac{1}{\sqrt{2\pi}}e^{-t^2/2}\,dt$$ denote the standard normal distribution function.

\begin{maintheorem}[Central limit theorem]\label{mainthm:clt}
Fix $g\ge1$. As $n\to\infty$, the variances of $X^H_{g,n}$ and $X^L_{g,n}$ satisfy
\[
 \Var(X^H_{g,n})=\nu n+\frac{(g-1)(11-4e)}{12(e-2)}+O(n^{-1/2}),
 \qquad
 \Var(X^L_{g,n})=\frac{n-g+1}{12}+O(n^{-1}),
\]
and, for $J\in\{H,L\}$,
\[
 \sup_{x\in\mathbb R}\left|\mathbb P\left(\frac{X^J_{g,n}-d_n/2}{\sqrt{\Var(X^J_{g,n})}}\le x\right)-\Phi(x)\right|=O(n^{-1/2}).
\]
The same conclusions hold for the full Betti numbers $b_j$, indexed by the cohomological degree $j$, with mean $d_n$ and with both variance expansions multiplied by $4$.
\end{maintheorem}

The leading variance coefficients are independent of the genus. The coefficient $\nu$ also governs the even Betti distributions of $\Mbar_{0,n}$ \cite{CKPNAS} and of the Fulton--MacPherson spaces \cite{CKTJM}, while $1/12$ comes from the Eulerian numbers. After normalization, the Betti polynomials have the asymptotic quasi-powers form $A(t)B(t)^n$, where $A(t)$ is the genus-dependent amplitude and $B(t)$ is the genus-independent base. The base determines the coefficient of $n$ in the variance, while the amplitude determines the constant term.

The Gaussian law of Theorem~\ref{mainthm:clt} concerns cumulative probabilities. We also prove a local limit theorem for the even Betti numbers (Theorem~\ref{thm:llt}), giving a Gaussian approximation that is uniform over all degrees and includes higher-order corrections. Together with the asymptotics of the even Betti totals, these local estimates yield asymptotic formulas for individual even Betti numbers. They involve the weighted count
\[
 c_g=\sum_{\Gamma}\frac1{|\Aut\Gamma|}
\]
of connected trivalent graphs $\Gamma$ of first Betti number $g$, with loops and multiple edges allowed. These numbers are determined by
\[
 \sum_{g\ge2}c_gx^{g-1}=\log\sum_{m\ge0}\frac{(6m)!}{2^{3m}(3m)!\,(3!)^{2m}(2m)!}\,x^m,
\]
so that $c_2=5/24$ and $c_3=5/16$ (Proposition~\ref{prop:core-constants}).

\begin{maintheorem}[Individual Betti numbers]\label{mainthm:betti}
Fix $g\ge1$ and let $x_{n,k}=k-d_n/2$. Uniformly over all integers $k$,
\[
 \begin{aligned}
 b_{2k}(\Mbar_{g,n})
 &=\frac{C_g(1)}{\sqrt{2\pi\nu}}\;n!\,(e-2)^{-n}\,n^{\frac{3g-6}{2}}
 \left[\exp\!\left(-\frac{x_{n,k}^2}{2\nu n}\right)+O(n^{-1/2})\right],\\
 b_{2k}(\Mbar_{g,0|n})
 &=\sqrt{\frac6\pi}\;A_g(1)\;n!\,n^{\frac{6g-9}{2}}
 \left[\exp\!\left(-\frac{6x_{n,k}^2}{n}\right)+O(n^{-1})\right],
 \end{aligned}
\]
where $A_1(1)=1/2$, $C_1(1)=1/4$, and, for $g\ge2$,
\[
 A_g(1)=\frac{c_g\,e^{-2(g-1)}}{\Gamma(3g-3)},\qquad
 C_g(1)=\frac{c_g\,e^{-(g-1)/2}}{\Gamma\!\left(\frac{3g-3}{2}\right)\,[2(e-2)]^{\frac{3g-3}{2}}}.
\]
\end{maintheorem}

The asymptotics of the even Betti totals are established in Theorem~\ref{thm:coefficients}. The odd cohomology is negligible in comparison: it is a proportion $O(n^{-1/2})$ of the cohomology of $\Mbar_{g,n}$ and $O(n^{-1})$ of that of $\Mbar_{g,0|n}$, and Theorem~\ref{mainthm:hodge} below identifies the exact orders. Consequently, the Euler characteristics are asymptotic to the even Betti totals, which answers a question of Kannan, Serpente and Yun \cite[Question~5.3]{KSY}.

\begin{maincorollary}[Euler characteristics]\label{maincor:euler}
For every fixed $g\ge2$,
\[
 \frac{\chi(\Mbar_{g,0|n})}{\chi(\Mbar_{g,n})}
 =\left(\frac{2(e-2)}{e}\right)^{\frac{3g-3}{2}}
 \frac{\Gamma\!\left(\frac{3g-3}{2}\right)}{\Gamma(3g-3)}\,
 n^{\frac{3g-3}{2}}(e-2)^n\bigl(1+O(n^{-1/2})\bigr).
\]
In particular, the ratio tends to zero.
\end{maincorollary}

For $g=1$, Kannan, Serpente and Yun proved that the ratio is asymptotic to $2(e-2)^n$ \cite[Corollary~5.2]{KSY}. Our methods recover this.

\subsection{Tautological classes}\label{subsec:intro-tautological}
Let $RH^{2k}(\Mbar_{g,n})\subset H^{2k}(\Mbar_{g,n})$ be the image of the tautological ring under the cycle class map. In genus zero, every class is tautological \cite{Keel}, and in genus one every even degree class is tautological \cite{Petersen}. For $g\ge2$ this fails: Graber and Pandharipande constructed a non-tautological algebraic class on $\Mbar_{2,20}$ \cite{GraberPandharipande}, and many further constructions are surveyed in \cite{FaberPandharipande}. Our next result shows that, nevertheless, non-tautological classes form a vanishing proportion of the even cohomology, of exact order $n^{-1}$.

\begin{maintheorem}[Tautological classes]\label{mainthm:taut}
Fix $g\ge2$. There are constants $0<c'_g\le C'_g$ such that, for all sufficiently large $n$,
\[
 \frac{c'_g}{n}\le
 1-\frac{\sum_k\dim RH^{2k}(\Mbar_{g,n})}{\sum_k b_{2k}(\Mbar_{g,n})}
 \le\frac{C'_g}{n}.
\]
Consequently, the distribution of the graded dimensions of $RH^*(\Mbar_{g,n})$ lies within total variation distance $O(n^{-1})$ of the even Betti distribution, and it satisfies the central limit theorem of Theorem~\ref{mainthm:clt} with the same mean $d_n/2$ and variance $\nu n+O(1)$.
\end{maintheorem}

The two bounds have different origins. For the upper bound, we show that every even class restricting to zero on the complement of the preimage of the boundary strata of dimension at most one in $\Mbar_g$ is tautological. We express such classes as sums of gluing pushforwards from products of genus-zero moduli spaces and at most one genus-one moduli space. The even cohomology of these products is tautological by the theorems of Keel and Petersen. The cohomology of that complement is a proportion $O(n^{-1})$ of the whole, by a version of the heavy/light substitution for compactly supported Betti numbers of open parts. For the lower bound, we exhibit classes of Hodge type $(p,q)$ with $p-q=22$. These come from stable cores with two genus-one vertices of valency one, called \emph{elliptic leaves}, each carrying the weight-twelve cusp form $\Delta$. Such classes cannot be tautological, since tautological classes have type $(k,k)$, and they form a proportion of exact order $n^{-1}$. In particular, for every $g\ge2$ and all large $n$, the even cohomology of $\Mbar_{g,n}$ contains non-tautological classes that are not even of Hodge type $(k,k)$.

\subsection{Hodge types off the diagonal}\label{subsec:intro-hodge}
For $Y=\Mbar_{g,n}$ or $Y=\Mbar_{g,0|n}$, choose a Hodge type $(p,q)$ with probability $h^{p,q}(Y)/\sum_j b_j(Y)$, and let $D^H_{g,n}$ and $D^L_{g,n}$ denote the Hodge difference $p-q$ in the heavy and light cases, respectively. For even $k\ge4$, let $\beta_k=\dim S_k(\mathrm{SL}_2(\mathbb Z))$ be the dimension of the space of cusp forms of weight $k$, and let $\mathbb I_m$ denote the modified Bessel function of the first kind. Let
\begin{equation}\label{eq:intro-Q}
 Z_Q=\sum_{\substack{m\ge3\\m\ \mathrm{odd}}}\beta_{m+1}\,m\,\mathbb I_m(2),\qquad q_m=q_{-m}=\frac{\beta_{m+1}\,m\,\mathbb I_m(2)}{2Z_Q}\quad(m\ge3\text{ odd}),
\end{equation}
and $q_d=0$ for all other integers $d$. These numbers define a probability distribution $Q$ on $\mathbb Z$.

\begin{maintheorem}[Hodge types off the diagonal]\label{mainthm:hodge}
Fix $g\ge2$. There are explicit constants $\gamma_g^H,\gamma_g^L>0$ such that
\[
 \mathbb P(D^H_{g,n}\ne0)=\gamma_g^Hn^{-1/2}+O(n^{-1}),\qquad \mathbb P(D^L_{g,n}\ne0)=\gamma_g^Ln^{-1}+O(n^{-2}).
\]
The same formulas hold for the probability that $p+q$ is odd, while the even-degree classes off the diagonal form a proportion $O(n^{-1})$ in the heavy case and $O(n^{-2})$ in the light case. For $J\in\{H,L\}$, the conditional distribution of $D^J_{g,n}$ given $D^J_{g,n}\ne0$ converges to $Q$. Moreover, for every $x\in\mathbb R$ and $d\in\mathbb Z$,
\[
 \mathbb P\!\left(\frac{p+q-d_n}{\sqrt{4\nu n}}\le x,\ D^H_{g,n}=d\ \middle|\ D^H_{g,n}\ne0\right)\longrightarrow\Phi(x)q_d,
\]
where $\Phi$ is the standard Gaussian distribution function. The same limit holds in the light case with $\sqrt{4\nu n}$ replaced by $\sqrt{n/3}$ and $D^H_{g,n}$ by $D^L_{g,n}$. Thus, conditionally on being off the Hodge diagonal, the normalized degree and the Hodge difference are asymptotically independent.
\end{maintheorem}

The off-diagonal classes therefore form a vanishing proportion of the cohomology, but their conditional Hodge difference has a universal limiting distribution, independent of the genus and common to the heavy and light spaces. Among the off-diagonal classes, the proportion of odd-degree classes tends to $1$. Since $\beta_{12}=1$ is the first nonzero cusp-form dimension, the limiting proportion with Hodge difference $\pm11$ is
\[
 q_{11}+q_{-11}=\frac{11\,\mathbb I_{11}(2)}{Z_Q}>0.9999.
\]

Geometrically, the leading off-diagonal classes come from boundary strata whose stable core is trivalent and rational except for a single elliptic leaf. For odd $m\ge3$, the local system $\operatorname{Sym}^{m-1}V$ on $\mathcal M_{1,1}$ contributes cusp forms of weight $m+1$ with Hodge difference $\pm m$ through Eichler--Shimura theory. Here $V$ is the rank-two local system given by the first cohomology of the universal elliptic curve. The remaining Tate factors shift $p$ and $q$ equally. This matches the special role of $\Mbar_{1,11}$ in the work of Canning, Larson and Payne, who proved that $H^{11}(\Mbar_{g,n})\ne0$ only for $g=1$ and $n\ge11$, and who used boundary pushforwards from $\Mbar_{1,11}$ to produce odd cohomology in higher genus \cite{CLP}. Remark~\ref{rem:skellam-interpretation} interprets the Bessel weights in $q_m$ through highest-weight multiplicities and Skellam probabilities.

\subsection{Ultra-log-concavity}\label{subsec:intro-ulc}
A positive sequence $(a_k)_{k=0}^d$ is \emph{log-concave} at an index $1\le k\le d-1$ if $a_k^2\ge a_{k-1}a_{k+1}$. For $r\ge0$, it is \emph{$r$-ultra-log-concave} at $k$ if the sequence $(a_k/\binom dk^r)_k$ is log-concave at $k$. For $r=0$ this is ordinary log-concavity, and for $r=1$ it is the classical ultra-log-concavity, studied for example through Lorentzian polynomials \cite{BrandenHuh}. Larger values of $r$ give stronger conditions. By hard Lefschetz, the even Betti numbers of $\Mbar_{g,n}$ and $\Mbar_{g,0|n}$ are unimodal. Near the middle degree we determine exactly how much stronger than log-concavity they are.

\begin{maintheorem}[Central ultra-log-concavity]\label{mainthm:ulc}
Fix $g\ge1$ and $0<c_0<\sqrt2$, and let
\[
 r_H=\frac1{4\nu}=\frac{3(e-2)}{2(3-e)}=3.82447\ldots,\qquad r_L=3.
\]
For $J\in\{H,L\}$, let $(a_k)$ be the even Betti sequence of $\Mbar_{g,n}$ if $J=H$ and of $\Mbar_{g,0|n}$ if $J=L$. Consider the central range of indices $k$ with $|k+j-d_n/2|\le c_0\sqrt{\Var(X^J_{g,n})\log n}$ for $j=-1,0,1$. For all sufficiently large $n$, the sequence $(a_k)$ is strictly $r$-ultra-log-concave at every such $k$ for each fixed $r\le r_J$, and fails to be $r$-ultra-log-concave at every such $k$ for each fixed $r>r_J$.
\end{maintheorem}

The central range has width of order $\sqrt{n\log n}$. The threshold has a simple probabilistic meaning. If a symmetric sequence of length about $n$ has variance about $vn$, comparison with the binomial coefficients $\binom nk$, whose normalized variance is about $n/4$, places the threshold at $1/(4v)$ (Lemma~\ref{lem:central-ulc}). The critical case $r=r_J$ is more delicate: there the terms of order $n^{-1}$ cancel, and the outcome depends on the fourth cumulant and on the constant terms of the mean and variance. It is not automatic. For the Losev--Manin spaces $\Mbar_{0,2|n}$, whose Betti numbers are the Eulerian numbers and whose threshold is also $3$, central $3$-ultra-log-concavity fails (Example~\ref{ex:losev-manin}).

\subsection{Method}\label{subsec:intro-method}
We work with the exponential generating functions $\mathscr L_g(w;u,v)$ and $\mathscr H_g(z;u,v)$ of the Hodge--Deligne polynomials of the light and heavy spaces. By \cite{KSY}, $\mathscr L_g(w)=\mathscr H_g(R(w,uv))$ for an explicit function $R$ built from genus zero. Writing $t=us$, let $T(z,t)$ denote the inverse of $R(w,t)$ in its first variable. This inverse is independent of the genus and has a square-root singularity at
\[
 \rho(t)=t^{\frac1{t-1}}-1-t^{-1},\qquad \rho(1)=e-2.
\]
This explains the factor $(e-2)^{-n}$ in the heavy asymptotics.

The light series are computed geometrically. Forgetting the markings of a light curve and stabilizing produces a stable core in $\Mbar_g$. The light curve is recovered from the core by inserting chains of rational curves at the nodes, which are parametrized by Losev--Manin spaces, and by placing markings on the smooth locus of the core. The chains contribute the Eulerian generating function which has a simple pole. Markings on the smooth core contribute entire functions; to see this, we bound the cohomology of fiber powers of universal curves exponentially, allowing nontrivial monodromy. A core with $|E|$ edges therefore produces a pole of order at most $|E|$. The maximal order $3g-3$ is attained exactly by the trivalent rational cores, which accounts for the constant $c_g$. In genus one there are no such cores, and cycles of rational curves produce a logarithm instead.

Singularity analysis then gives uniform quasi-powers estimates, and Hwang's theorem \cite{HwangQP} yields Theorem~\ref{mainthm:clt}. For the local limit theorem we refine Hwang's Fourier method \cite{HwangLLT}. Poincar\'e duality makes the distributions symmetric about their means, so all odd cumulants of the centered variables vanish. The resulting local expansion is precise enough to compare neighboring Betti numbers at the critical exponent.

For Hodge types we substitute $(u,v)=(-s\xi,-s\xi^{-1})$, so that $\xi$ records $p-q$. Rational cores are Hodge--Tate, and the dimension of a boundary stratum bounds the pole order of its contribution. The first non-Tate contributions therefore come from the strata of dimension one with an elliptic leaf. Their leading pole coefficients involve the cusp forms of weight $m+1$, through the computation of the Hodge--Deligne polynomials in genus one \cite{Kiem}, and the Bessel functions arise from summing over the markings on the leaf. The same analysis, applied to two elliptic leaves, gives the lower bound in Theorem~\ref{mainthm:taut}.

\subsection{Layout}
Section~\ref{sec:preliminaries} collects the Hodge-theoretic and analytic tools. Section~\ref{sec:light} computes the singular structure of the Hodge generating series. Section~\ref{sec:clt} proves the coefficient asymptotics, Theorem~\ref{mainthm:clt} and Corollary~\ref{maincor:euler}. Section~\ref{sec:local} proves the local limit theorem, Theorem~\ref{mainthm:betti} and Theorem~\ref{mainthm:ulc}. Section~\ref{sec:hodge-refined} develops the Hodge-refined series and the analysis of strata by dimension, and proves Theorem~\ref{mainthm:hodge}. Section~\ref{sec:tautological} proves Theorem~\ref{mainthm:taut}.

The main theorems correspond to results in the body as follows: Theorem~\ref{mainthm:clt} to Theorem~\ref{thm:clt} and Corollary~\ref{cor:full-degree}; Theorem~\ref{mainthm:betti} to Corollary~\ref{cor:even-betti-estimates}; Corollary~\ref{maincor:euler} to Corollary~\ref{cor:euler-ratio}; Theorem~\ref{mainthm:taut} to Theorem~\ref{thm:tautological-dominance}; Theorem~\ref{mainthm:hodge} to Theorem~\ref{thm:hodge-asymmetry}, where $\gamma^H_g$ and $\gamma^L_g$ are given in \eqref{eq:hodge-rarity-constants}; and Theorem~\ref{mainthm:ulc} to Theorems~\ref{thm:heavy-ulc} and~\ref{thm:light-ulc}.

\section{Preliminaries}\label{sec:preliminaries}
This section collects the Hodge-theoretic identities behind our generating series and the analytic tools used to extract asymptotics from them.

\subsection{Hodge--Deligne polynomials and finite quotients}\label{subsec:hodge-polynomials}
For a complex algebraic variety, or a Deligne--Mumford stack, $Y$, let $E(Y;u,v)$ denote the Hodge--Deligne polynomial of its compactly supported rational cohomology. For an automorphism $\alpha$ of finite order, the equivariant version is
\[
 E_\alpha(Y;u,v)=\sum_{i,p,q}(-1)^i
 \Tr\left(\alpha\mid
 \operatorname{Gr}_F^p\operatorname{Gr}_{p+q}^W H_c^i(Y;\mathbb C)
 \right)u^pv^q.
\]
If a finite group $G$ acts on $Y$, the rational cohomology of the quotient is the invariant part, compatibly with mixed Hodge structures, and averaging traces gives
\begin{equation}\label{eq:trace-average}
 E(Y/G;u,v)=\frac1{|G|}\sum_{\alpha\in G}E_\alpha(Y;u,v).
\end{equation}

The same argument applies to stacks, because the higher cohomology of a finite group with rational coefficients vanishes. The polynomials $E_\alpha$ are additive for invariant stratifications, by the long exact sequence of a closed invariant subspace and its complement, and multiplicative for products, by the K\"unneth formula. If $\alpha$ permutes a family of strata without fixing any of them, the trace of $\alpha$ on their union is zero.

For either $Y=\Mbar_{g,n}$ or $Y=\Mbar_{g,0|n}$, smoothness and properness of $Y$ imply that $H^m(Y;\mathbb Q)$ carries a pure Hodge structure of weight $m$. Hence the coefficient of $u^pv^q$ in $E(Y;u,v)$ equals $(-1)^{p+q}h^{p,q}(Y)$, and
\begin{equation}\label{eq:parity-extraction}
 P_Y^{\ev}(t)=\frac{E(Y;\sqrt t,\sqrt t)
                         +E(Y;-\sqrt t,-\sqrt t)}2,
 \qquad P_Y(q)=E(Y;-q,-q),
\end{equation}
where $P_Y(q)=\sum_j b_j(Y)q^j$ is the Poincar\'e polynomial and $P_Y^{\ev}(t)=\sum_k b_{2k}(Y)t^k$ is the even degree Poincar\'e polynomial. The average in \eqref{eq:parity-extraction} retains all Hodge types with $p+q$ even and does not depend on the choice of square root.

\subsection{Singularity analysis}\label{subsec:singularity-analysis}
We write $[z^n]f(z)$ for the coefficient of $z^n$ in a power series $f$. All our generating series are exponential in the number of markings,
\[
 \mathcal A(z,t)=\sum_{n\ge0}A_n(t)\frac{z^n}{n!},
 \qquad A_n(t)=n![z^n]\mathcal A(z,t),
\]
so that a product of two such series distributes the labels between the factors, with the binomial coefficients supplied automatically.

A \emph{dominant singularity} of a power series is a singularity on its circle of convergence. For the logarithmic and algebraic singularities considered here, a unique dominant singularity determines the coefficient asymptotics: its location gives the exponential growth rate, and the local form of the function gives the polynomial factor. The basic transfer models are
\begin{equation}\label{eq:transfer-models}
 [z^n]\log\frac1{1-z/\rho}=\frac{\rho^{-n}}n,\qquad
 [z^n](1-z/\rho)^{-\alpha}
   \sim\frac{\rho^{-n}n^{\alpha-1}}{\Gamma(\alpha)},
\end{equation}
valid for $n\ge1$ and for $\alpha\notin\{0,-1,-2,\ldots\}$, respectively. Thus a pole of order $m$ contributes $n^{m-1}$, and a square-root term contributes $n^{-3/2}$.

To transfer a local expansion together with its error term, one needs analytic continuation slightly beyond the circle of convergence. A \emph{$\Delta$-domain} at $z=\rho$ is a set of the form
\[
 |z/\rho|<1+\eta,\qquad |\arg(z/\rho-1)|>\theta,
 \qquad \eta>0,\quad 0<\theta<\pi/2.
\]
The transfer theorem \cite[Theorem~VI.4]{FS} converts a local expansion with remainder on such a domain into a coefficient expansion with remainder. We use it with a complex parameter $t$ near $1$. If the function depends holomorphically on $t$, the domains are obtained from one fixed $\Delta$-domain by rescaling with $\rho(t)$, and both the local remainder and the bounds away from the singularity are uniform in $t$, then the coefficient estimates are uniform in $t$ \cite[Lemma~IX.2]{FS}.

\subsection{Probability generating functions and quasi-powers}\label{subsec:pgf-quasi-powers}
Let $A(t)=\sum_{k=0}^d a_k t^k$ be a polynomial with nonnegative coefficients and $A(1)>0$. It defines an integer-valued random variable $Y$ with $\mathbb P(Y=k)=a_k/A(1)$, whose \emph{probability generating function} (PGF) is
\begin{equation}\label{eq:prelim-pgf}
 G_Y(t)=\mathbb E t^Y=\frac{A(t)}{A(1)}.
\end{equation}
Applied to an even Betti polynomial, this gives the distribution of half the cohomological degree. For a function $f$ holomorphic near $1$ with $f(1)\ne0$, set
\begin{equation}\label{eq:analytic-moments}
 \mathfrak m(f)=\frac{f'(1)}{f(1)},\qquad
 \mathfrak v(f)=\frac{f''(1)}{f(1)}+\frac{f'(1)}{f(1)}
                   -\left(\frac{f'(1)}{f(1)}\right)^2.
\end{equation}
For a local holomorphic logarithm of $f$, the chain rule gives
\begin{equation}\label{eq:analytic-moments-chain}
 \left.\frac{d}{ds}\log f(e^s)\right|_{s=0}=\mathfrak m(f),
 \qquad
 \left.\frac{d^2}{ds^2}\log f(e^s)\right|_{s=0}=\mathfrak v(f),
\end{equation}
so $\mathbb E Y=\mathfrak m(G_Y)$ and $\Var(Y)=\mathfrak v(G_Y)$ \cite[p.~645]{FS}. The expressions $ \mathfrak m(f)$ and $\mathfrak v (f)$ are unchanged when $f$ is multiplied by a nonzero constant. If $A$ is palindromic, $A(t)=t^dA(t^{-1})$, then $Y$ and $d-Y$ have the same law, so $\mathbb E Y=d/2$ exactly. By Poincar\'e duality this applies to the even Betti polynomials of $\Mbar_{g,n}$ and $\Mbar_{g,0|n}$, which therefore have exact mean $d_n/2$.

Our coefficient estimates, divided by their values at $t=1$, have the following form. We use Hwang's theorem \cite{HwangQP} as stated in \cite[Theorem~IX.8]{FS}.

\begin{theorem}[Quasi-powers theorem]\label{thm:quasi-powers}
Let $Y_n$ be nonnegative integer-valued random variables whose PGFs $p_n(t)$ are holomorphic on a fixed complex neighborhood of $t=1$. Suppose that, uniformly on this neighborhood,
\begin{equation}\label{eq:quasi-powers-pgf}
 p_n(t)=A(t)B(t)^{\beta_n}\bigl(1+O(\varepsilon_n)\bigr),
 \qquad \beta_n\to\infty,\quad \varepsilon_n\to0,
\end{equation}
where $A$ and $B$ are holomorphic and independent of $n$, $A(1)=B(1)=1$, $\mathfrak v(B)>0$, and the power is defined with the branch of $\log B$ vanishing at $1$. Then the mean $\mu_n$ and variance $\sigma_n^2$ of $Y_n$ satisfy
\[
 \mu_n=\beta_n\mathfrak m(B)+\mathfrak m(A)+O(\varepsilon_n),\qquad
 \sigma_n^2=\beta_n\mathfrak v(B)+\mathfrak v(A)+O(\varepsilon_n),
\]
and
\[
 \sup_{x\in\mathbb R}\left|
 \mathbb P\left(\frac{Y_n-\mu_n}{\sigma_n}\le x\right)
 -\Phi(x)\right|=O(\beta_n^{-1/2}+\varepsilon_n),
\]
where $\Phi$ is the standard normal distribution function.
\end{theorem}

The factors $A$ and $B$ need not be PGFs themselves. The coefficient $\mathfrak v(B)$ of the leading term of the variance is determined by the base $B$ alone.

\subsection{Cumulants and Fourier inversion}\label{subsec:cumulants}
To go beyond the mean and the variance we use cumulants, in the logarithmic coordinate $t=e^s$. For an integer-valued random variable $Y$ with finite support, the \emph{moment generating function} $M_Y(s)=\mathbb E e^{sY}=G_Y(e^s)$ is entire and satisfies $M_Y(0)=1$. On a neighborhood of $0$, its logarithm
\[
 \Lambda_Y(s)=\log M_Y(s)=\sum_{j\ge1}\kappa_j(Y)\frac{s^j}{j!},\qquad \Lambda_Y(0)=0,
\]
is the \emph{cumulant generating function}, and $\kappa_j(Y)$ is the $j$th \emph{cumulant}. If $\mu=\mathbb E Y$ and $\sigma^2=\Var(Y)$, then
\begin{equation}\label{eq:prelim-cumulants}
 \kappa_1(Y)=\mu,\qquad \kappa_2(Y)=\sigma^2,\qquad
 \kappa_3(Y)=\mathbb E(Y-\mu)^3,\qquad
 \kappa_4(Y)=\mathbb E(Y-\mu)^4-3\sigma^4.
\end{equation}
Cumulants are additive for independent sums, and the cumulants of order at least three vanish for a Gaussian variable. They therefore measure the departure from a Gaussian law. If $G_Y$ is a palindromic polynomial of degree $d$, then $M_Y(s)=e^{ds}M_Y(-s)$. The centered cumulant generating function $\Lambda_Y(s)-ds/2$ is then even, and all odd cumulants of $Y-d/2$ vanish. This is the probabilistic form of Poincar\'e duality that we exploit in Section~\ref{sec:local}.

The quasi-powers hypothesis controls all cumulants. Suppose \eqref{eq:quasi-powers-pgf} holds with $\beta_n=n$, and let
\[
 \lambda(s)=\log B(e^s),\qquad \psi(s)=\log A(e^s),
\]
with both logarithms vanishing at $s=0$. Then \eqref{eq:quasi-powers-pgf} becomes
\begin{equation}\label{eq:quasi-powers}
 M_n(s)=\exp\bigl(n\lambda(s)+\psi(s)\bigr)
                  \bigl(1+O(\varepsilon_n)\bigr)
\end{equation}
uniformly on a fixed complex neighborhood of $0$, where $M_n=M_{Y_n}$ and $\lambda''(0)=\mathfrak v(B)>0$. For large $n$, the function $M_n(s)\exp(-n\lambda(s)-\psi(s))$ is uniformly close to $1$, so it has a holomorphic logarithm $R_n(s)=O(\varepsilon_n)$ with $R_n(0)=0$. Cauchy's estimates on a smaller disc give $R_n^{(j)}(s)=O_j(\varepsilon_n)$ for every fixed $j$, and hence
\begin{equation}\label{eq:cumulants}
 \kappa_j(Y_n)=n\lambda^{(j)}(0)+\psi^{(j)}(0)+O_j(\varepsilon_n).
\end{equation}
Here $O_j(\varepsilon_n)$ denotes a quantity bounded by $C_j\varepsilon_n$, where $C_j$ may depend on $j$ but not on $n$. In particular, every fixed cumulant is $O(n)$, so the standardized fourth cumulant $\kappa_4(Y_n)/\Var(Y_n)^2$ is $O(n^{-1})$. For symmetric laws this is the first correction to the Gaussian approximation.

A central limit theorem controls cumulative probabilities, whereas a local limit theorem controls the individual probabilities $\mathbb P(Y_n=k)$. In our setting, the error $O(n^{-1/2})$ in the central limit theorem is of the same order as the central probabilities themselves, so a separate local argument is needed. Lattice effects can prevent a local limit theorem: if $B_N$ is binomial with parameters $N$ and $1/2$, then $2B_N$ satisfies a central limit theorem but has zero probability at every odd integers.

For an integer-valued variable $Y$, Fourier inversion, or Cauchy's coefficient formula on the unit circle, gives
\begin{equation}\label{eq:prelim-fourier}
 \mathbb P(Y=k)=\frac1{2\pi}\int_{-\pi}^{\pi}
                       e^{-ik\theta}M_Y(i\theta)\,d\theta.
\end{equation}
In Section~\ref{sec:local} we evaluate this integral using the cumulant expansion near $\theta=0$ and an exponential bound on $|M_Y(i\theta)|$ away from $\theta=0$. The latter bound rules out lattice effects.
\section{Hodge series and the heavy/light substitution}\label{sec:light}

\subsection{Weighted pointed curves and generating series}
We recall Hassett's moduli spaces of weighted pointed stable curves \cite{Hassett}. Let $\mathcal{A}=(a_1,\ldots,a_n)$ be a vector of rational weights with $0<a_i\le1$ and $2g-2+\sum_i a_i>0$. A connected projective nodal curve $C$ of arithmetic genus $g$ with labeled markings $p_1,\ldots,p_n$ in its smooth locus is \emph{$\mathcal{A}$-stable} if
\begin{enumerate}
\item whenever the markings indexed by $I\subset \{1,\ldots,n\}$ coincide, $\sum_{i\in I}a_i\le1$, and
\item the $\mathbb Q$-divisor $K_C+\sum_i a_i p_i$ is ample.
\end{enumerate}
The second condition means that
\[
 2g_v-2+d_v+\sum_{p_i\in C_v}a_i>0
\]
for every irreducible component $C_v$, where $g_v$ is the genus of its normalization and $d_v$ is the number of branches of nodes on $C_v$. A self-node contributes two to $d_v$, and coincident markings are counted separately. These curves are parametrized by a smooth proper Deligne--Mumford stack $\Mbar_{g,\mathcal{A}}$. For $g\ge1$ and $n\ge1$ we use the two extreme weight vectors
\[
 \Mbar_{g,n}=\Mbar_{g,(1,\ldots,1)},\qquad
 \Mbar_{g,0|n}=\Mbar_{g,(1/n,\ldots,1/n)}.
\]
In the heavy space $\Mbar_{g,n}$ the markings are distinct and Hassett stability is ordinary stability. In the light space $\Mbar_{g,0|n}$ the total weight is one, so any set of markings may coincide at a smooth point. Both spaces have dimension $d_n=n+3g-3$. For $g\ge2$ and $n=0$, both are understood to be $\Mbar_g$. For $g=1$, the unmarked space is unstable, and we set the constant terms of all generating series below equal to zero.

Define the exponential generating series of Hodge--Deligne polynomials
\begin{equation}\label{eq:hodge-series}
 \mathscr L_g(w;u,v)=\sum_{n\ge0}
 E(\Mbar_{g,0|n};u,v)\frac{w^n}{n!},\qquad
 \mathscr H_g(z;u,v)=\sum_{n\ge0}
 E(\Mbar_{g,n};u,v)\frac{z^n}{n!},
\end{equation}
the even Betti polynomials
\begin{equation}\label{eq:betti-polynomials}
 L_{g,n}(t)=\sum_{k=0}^{d_n}b_{2k}(\Mbar_{g,0|n})t^k,\qquad
 H_{g,n}(t)=\sum_{k=0}^{d_n}b_{2k}(\Mbar_{g,n})t^k,
\end{equation}
and their exponential generating series
\[
 \mathcal L_g(w,t)=\sum_{n\ge0}L_{g,n}(t)\frac{w^n}{n!},\qquad
 \mathcal H_g(z,t)=\sum_{n\ge0}H_{g,n}(t)\frac{z^n}{n!}.
\]
Script letters $\mathscr L,\mathscr H$ denote full Hodge series, and calligraphic letters $\mathcal L,\mathcal H$ their even Betti parts.

\subsection{The heavy/light substitution}
The heavy and light series are related by a substitution that is independent of the genus and comes from genus zero. The spaces $\Mbar_{0,n+1}$ have no odd cohomology, and $H^{2k}(\Mbar_{0,n+1})$ is of Hodge type $(k,k)$. Let
\[
 \varphi(z,t)=z+\sum_{n\ge2}\frac{z^n}{n!}
                    \sum_{k=0}^{n-2}b_{2k}(\Mbar_{0,n+1})t^k.
\]
Summing over stable trees gives the functional equation \cite[p.~194, Theorem~4.2.1]{Manin}
\begin{equation}\label{eq:phi-functional}
 (1+\varphi)^t=t^2\varphi-t(t-1)z+1.
\end{equation}
Equivalently, by the non-equivariant specialization of Getzler's Legendre-transform relation \cite[Theorem~5.9]{Getzler}, the formal inverse of $\varphi$ in its first variable is
\[
 \varphi^{\langle-1\rangle}(x,t)
 =x-\sum_{n\ge2}E(\mathcal M_{0,n+1};u,v)\frac{x^n}{n!}
 =\frac{t^2x-(1+x)^t+1}{t(t-1)},\qquad t=uv.
\]

Kannan, Serpente and Yun proved that the light Hodge series is obtained by substituting a universal function into the heavy one \cite[Corollary~B]{KSY}. In the all-light case it reads
\begin{equation}\label{eq:KSYsub}
 \mathscr L_g(w;u,v)=\mathscr H_g(R(w,uv);u,v).
\end{equation}
The substitution combines two geometric steps \cite[Propositions~4.6--4.7]{KSY}. First, a heavy curve is obtained from a curve with distinct markings and no rational tails by attaching genus-zero tails at the markings. Since $\varphi$ records these tails, substituting $\varphi^{\langle-1\rangle}(x,uv)$ for $z$ in $\mathscr H_g$ gives the series of the intermediate locus, where $x$ counts the distinct markings. Second, allowing markings to coincide replaces each distinct marking by a nonempty set of light labels, whose exponential generating series is $e^w-1$. The combined substitution is therefore
\begin{equation}\label{eq:inverse-definition}
 R(w,t)=\varphi^{\langle-1\rangle}(e^w-1,t)
       =\frac{t^2e^w-t^2-e^{tw}+1}{t(t-1)},\qquad
 R(w,1)=(2-w)e^w-2.
\end{equation}
Since $R(0,t)=0$ and $R_w(0,t)=1$, $R$ has a unique local inverse $T(z,t)$ in its first variable satisfying $T(0,t)=0$. Inverting the composition in \eqref{eq:inverse-definition} gives
\begin{equation}\label{eq:genus-zero-log}
 T(z,t)=\log(1+\varphi(z,t)),
\end{equation}
with the branch of the logarithm vanishing at $1$. The coefficients of $R$ as a series in $w$ lie in $\mathbb Q[t]$, so formal inversion shows that the coefficients of $T$ as a series in $z$ are polynomials in $t$, and $T$ is holomorphic near $z=0$, uniformly for $t$ near $1$.

\begin{proposition}\label{prop:substitution}
For every $g\ge1$,
\begin{equation}\label{eq:substitution}
 \mathscr H_g(z;u,v)=\mathscr L_g(T(z,uv);u,v),
 \qquad \mathcal H_g(z,t)=\mathcal L_g(T(z,t),t).
\end{equation}
\end{proposition}

\begin{proof}
Substituting $w=T(z,uv)$ into the \eqref{eq:KSYsub} and using $R(T(z,uv),uv)=z$ gives the first identity. The two specializations $(u,v)=(\sqrt t,\sqrt t)$ and $(-\sqrt t,-\sqrt t)$ in \eqref{eq:parity-extraction} have the same product $uv=t$, hence the same substituted argument $T(z,t)$. Averaging them gives the second identity.
\end{proof}

The singularity of $T$ governs all heavy asymptotics. Let
\begin{equation}\label{eq:universal-constants}
 \tau(t)=\frac{\log t}{t-1},\qquad
 \rho(t)=t^{\frac{1}{t-1}}-1-t^{-1},
\end{equation}
where $\log t$ vanishes at $t=1$, and $\tau(1)=1$, $\rho(1)=e-2$ by continuity.

\begin{lemma}\label{lem:inverse}
For complex $t$ sufficiently close to $1$, the function $T(z,t)$ has a unique dominant singularity at $z=\rho(t)$, and it continues analytically to $\Delta_t=\rho(t)\Delta_0$, where $\Delta_0$ is a fixed $\Delta$-domain at $1$. As $z\to\rho(t)$ in $\Delta_t$,
\begin{equation}\label{eq:T-expansion}
 T(z,t)=\tau(t)-(2\rho(t))^{\frac12}\,t^{-\frac{t}{2(t-1)}}
 \left(1-\frac{z}{\rho(t)}\right)^{\frac12}
 +O\!\left(1-\frac{z}{\rho(t)}\right),
\end{equation}
uniformly for $t$ in a fixed complex neighborhood of $1$. Here $(2\rho(t))^{1/2}$ is positive at $t=1$, where $t^{-\frac{t}{2(t-1)}}$ takes the value $e^{-1/2}$.
\end{lemma}

\begin{proof}
A direct computation gives
\[
 R_w(\tau(t),t)=0,\qquad R(\tau(t),t)=\rho(t),
\]
so $w=\tau(t)$ is a critical point of $R(\cdot,t)$ with critical value $\rho(t)$. As in \cite[Equations~(3.3)--(3.8)]{CKPNAS}, the inverse function theorem at this quadratic critical point gives a convergent expansion
\[
 1+\varphi(z,t)=t^{\frac{1}{t-1}}
 -\left(2\rho(t)t^{\frac{2-t}{t-1}}\right)^{\frac12}
  \left(1-\frac{z}{\rho(t)}\right)^{\frac12}
 +O\!\left(1-\frac{z}{\rho(t)}\right),
\]
whose coefficients depend holomorphically on $t$ near $1$. Taking logarithms, with value $\tau(t)$ at $t^{1/(t-1)}\ne0$, gives \eqref{eq:T-expansion}.

It remains to show that $\rho(t)$ is the only dominant singularity. At $t=1$, the identity $R(w,1)=w-\sum_{n\ge3}\frac{n-2}{n!}w^n$ shows recursively that $T(z,1)$ has nonnegative Taylor coefficients. Since $R_w(w,1)=(1-w)e^w>0$ on $[0,1)$, $R(0,1)=0$ and $R(1,1)=e-2$, the real inverse is analytic on $[0,e-2)$ and tends to $1$ as $z$ approaches $e-2$ from below. By Pringsheim's theorem \cite[Theorem~IV.6]{FS}, the radius of convergence is $e-2$. The coefficients sum to $T(e-2,1)=1$, so the series converges absolutely on $|z|\le e-2$. If $|z|=e-2$ and $R_w(T(z,1),1)=0$, then $T(z,1)=1$ and $z=R(1,1)=e-2$. Hence the inverse function theorem continues $T(\cdot,1)$ across every other boundary point. Covering the circle $|z|=e-2$ outside a neighborhood of $e-2$ by finitely many such inverse-function neighborhoods and applying the implicit function theorem in the parameter $t$, we obtain the continuation to $\Delta_t=\rho(t)\Delta_0$ for all $t$ sufficiently close to $1$. The uniformity of the error term follows from the holomorphic dependence of the local expansion on $t$.
\end{proof}

The square-root singularity of $T$ arises from the critical point $w=\tau(t)$ of $R$. To make the vanishing factor of $R_w$ explicit, we differentiate \eqref{eq:inverse-definition} and write
\[
 R_w(w,t)=\frac{te^w-e^{tw}}{t-1}
 =e^w\left(1-\frac{e^{(t-1)w}-1}{t-1}\right).
\]
We write the second factor as $1-W(w,t)$, where
\begin{equation}\label{eq:W}
 W(w,t)=\frac{e^{(t-1)w}-1}{t-1},\qquad W(w,1)=w.
\end{equation}
The function $W$ is entire and will enter the singularity analysis of the light series. For $t$ near $1$, the pole of $1/(1-W(w,t))$ closest to the origin is $w=\tau(t)$. Indeed $W(\tau(t),t)=1$ and $W_w(\tau(t),t)=t\ne0$. Thus this pole occurs at the critical point of $R$ whose image $\rho(t)$ is the dominant singularity of its inverse $T$. Note that $1/(1-W)$ is the exponential generating function of the Eulerian polynomials \cite[Example~III.25]{FS}.

\subsection{The logarithmic singularity in genus one}\label{subsec:light-genus-one}
Following \cite[Section~4]{Kiem}, we analyze the genus-one light spaces through the morphisms $\Mbar_{1,0|n}\to\Mbar_{1,1}$, $n\ge1$, that forget all markings but the first and stabilize.

\begin{lemma}\label{lem:light-genus-one}
The full Hodge series has the form
\begin{equation}\label{eq:genus-one-hodge-series}
 \mathscr L_1(w;u,v)
 =-\frac12\log(1-W(w,uv))+\mathscr J_1(w;u,v),
\end{equation}
where the logarithm vanishes at $w=0$ and $\mathscr J_1$ is entire in $w$, locally uniformly in $(u,v)$, and holomorphic in all variables.
\end{lemma}

\begin{proof}
We decompose $\Mbar_{1,0|n}$ into the fiber over the boundary point of $\Mbar_{1,1}$ and the preimage of $\mathcal M_{1,1}$.

\emph{The boundary fiber.} Its curves have only rational components. Since the total weight is one, stability forces every component to have at least two node branches. The genus condition then forces the dual graph to be a cycle, and stability requires each component to carry a marking. For a cycle of $j\ge3$ components, the Hodge contribution is
\[
 \left\{\begin{matrix}n\\j\end{matrix}\right\}
       \frac{(j-1)!}{2}(uv-1)^{n-j}
\]
by \cite[Proposition~4.3]{Kiem}. Here the Stirling number counts the partitions of the $n$ labels into $j$ blocks, $(j-1)!/2$ counts their cyclic orders up to reflection, and a component carrying $k$ labels, with its two attaching points at $0$ and $\infty$, contributes $(\mathbb C^*)^k/\mathbb C^*\cong(\mathbb C^*)^{k-1}$, with Hodge polynomial $(uv-1)^{k-1}$. Using $\sum_{n}\left\{\begin{smallmatrix}n\\j\end{smallmatrix}\right\}x^n/n!=(e^x-1)^j/j!$ \cite[Appendix~A.8]{FS} and writing $W=W(w,uv)$, the cycles with at least three components contribute
\[
 \sum_{j\ge3}\sum_{n\ge j}
 \left\{\begin{matrix}n\\j\end{matrix}\right\}
 \frac{(j-1)!}{2}(uv-1)^{n-j}\frac{w^n}{n!}
 =\frac12\sum_{j\ge3}\frac{W^j}{j}=-\frac12\log(1-W)-\frac W2-\frac{W^2}{4}.
\]
By \cite[Proposition~4.3]{Kiem}, the cycles with one and two components contribute
\[
 \frac{(uv+1)^{n-1}+(uv-1)^{n-1}}{2}\quad(n\ge1),\qquad
 (2^{n-1}-1)\frac{(uv+1)^{n-2}+(uv-1)^{n-2}}{2}\quad(n\ge2).
\]
If $|uv\pm1|\le A$ with $A\ge1$, these are bounded by $(2A)^n$, so their exponential generating series are dominated on $|w|\le R$ by $e^{2AR}$ and are entire in $w$ and holomorphic in $(w,u,v)$.

\emph{The interior.} Let $I_n(u,v)$ be the Hodge polynomial of the preimage of $\mathcal M_{1,1}$. By \cite[Proposition~4.1]{Kiem}, $I_n(u,v)=I_n^{\mathrm{even}}(uv)+I_n^{\mathrm{odd}}(u,v)$, where the superscripts refer to the parity of the total Hodge degree and
\[
 \begin{aligned}
 I_n^{\mathrm{even}}(t)&=\sum_{j=0}^{\lfloor(n-1)/2\rfloor}
 \binom{n-1}{2j}(1+t)^{n-1-2j}
 \left(\alpha_{jj}t^{j+1}-\sum_{i=0}^{j-1}\alpha_{ij}t^i\right),\\
 I_n^{\mathrm{odd}}(u,v)&=-\sum_{j=0}^{\lfloor(n-1)/2\rfloor}
 \binom{n-1}{2j}(1+uv)^{n-1-2j}\sum_{i=0}^{j-1}\alpha_{ij}\beta_{2j-2i+2}(uv)^i
 \bigl(u^{2j-2i+1}+v^{2j-2i+1}\bigr),
 \end{aligned}
\]
with
\[
 \alpha_{ij}=\binom{2j}{i}-\binom{2j}{i-1},\qquad
 \beta_\ell=\begin{cases}
 \lfloor\ell/12\rfloor-1,&\ell\equiv2\pmod{12},\\
 \lfloor\ell/12\rfloor,&\text{otherwise}.
 \end{cases}
\]
For even $\ell\ge4$, $\beta_\ell=\dim S_\ell(\mathrm{SL}_2(\mathbb Z))$. 
Put $t=uv$ and $M=\max(1,|u|,|v|)$. Since $0\le\alpha_{ij}\le4^j$ for $0\le i\le j$, $0\le\beta_{2j-2i+2}\le n$ and $2j\le n$, we have
\[
 \left|\alpha_{jj}t^{j+1}-\sum_{i=0}^{j-1}\alpha_{ij}t^i\right|
 \le(j+1)4^jM^{2j+2}\le nM^2(2M)^{2j},
\]
and the inner sum in $I_n^{\mathrm{odd}}$ is at most $n^2M(2M)^{2j}$ in absolute value. As $|1+uv|\le1+M^2$, the binomial theorem gives
\[
 |I_n(u,v)|\le K n^2M^2\sum_{j}
 \binom{n-1}{2j}(1+M^2)^{n-1-2j}(2M)^{2j}
 \le K n^2M^2(M+1)^{2n-2}
\]
for an absolute constant $K$. Hence the exponential generating series of $I_n^{\mathrm{even}}(uv)$ and $I_n^{\mathrm{odd}}(u,v)$ converge uniformly on compact subsets of $(w,u,v)$. Collecting the interior series, the one- and two-component cycles and $-W/2-W^2/4$ gives $\mathscr J_1$.
\end{proof}

Averaging the specializations $(u,v)=(\sqrt t,\sqrt t)$ and $(-\sqrt t,-\sqrt t)$ cancels $I_n^{\mathrm{odd}}$, while all other contributions depend only on $uv=t$. Hence
\begin{equation}\label{eq:genus-one-series}
 \mathcal L_1(w,t)=-\frac12\log(1-W(w,t))+\mathcal{J}_1(w,t),
\end{equation}
with
\[
 \mathcal{J}_1(w,t)=\sum_{n\ge1}I_n^{\mathrm{even}}(t)\frac{w^n}{n!}
 +\frac{e^{(t+1)w}-1}{2(t+1)}
 +\frac{(e^{(t+1)w}-1)^2}{4(t+1)^2}.
\]
The last two terms come from the cycles with one and two components after subtraction of $W/2$ and $W^2/4$, and they are extended by continuity to $t=-1$. The estimates above show that $\mathcal{J}_1$ is entire in $w$, locally uniformly in $t$.

Since the even part of $E(\Mbar_{1,0|n};u,v)$ depends only on $uv$, the group $H^{2k}(\Mbar_{1,0|n})$ is of Hodge type $(k,k)$. The coefficients of $T(z,uv)$ are polynomials in $uv$, so Proposition~\ref{prop:substitution} shows that $H^{2k}(\Mbar_{1,n})$ is of Hodge type $(k,k)$ as well.

\subsection{Stable cores in higher genus}\label{subsec:light-higher-genus}
Now let $g\ge2$. We describe a light curve through its stable core, the markings on the smooth part of the core, and chains of rational curves inserted at the nodes, and we translate this description into generating functions.

Forgetting the markings and stabilizing defines $\pi_n:\Mbar_{g,0|n}\to\Mbar_g$. The dual graph $\Gamma$ of the stabilized curve $C$ has no legs, and its vertex genera $g_v$ and valencies $d_v$ satisfy $2g_v-2+d_v>0$ and $g=\sum_v g_v+b_1(\Gamma)$. After labeling the half-edges (flags), the cores of type $\Gamma$ are parametrized by
\[
 S_\Gamma=\prod_{v\in V(\Gamma)}\mathcal M_{g_v,d_v}:
\]
a point of $S_\Gamma$ is a collection of smooth pointed curves, and gluing the two flag points of each edge produces $C$. The stratum of $\Mbar_g$ with dual graph $\Gamma$ is $\mathcal M_\Gamma=[S_\Gamma/\Aut\Gamma]$, and
\begin{equation}\label{eq:stratum-dimension}
 \dim\mathcal M_\Gamma=\sum_v(3g_v-3+d_v)=3g-3-|E(\Gamma)|.
\end{equation}

A light curve with core $C$ is obtained by distributing the $n$ labels between the smooth locus of $C$ and its nodes. Markings on the smooth locus may coincide. Stability forces the rational components contracted by $\pi_n$ to form chains at the nodes, each component carrying at least one marking. At a node assigned $k\ge1$ labels we therefore insert a chain parametrized by the Losev--Manin space $\mathrm{LM}_k=\Mbar_{0,2|k}$ \cite{LosevManin}, whose two heavy markings are the ordered endpoints and whose light markings carry the assigned labels. For $k=0$ the node is left unchanged, and we set $\mathrm{LM}_0=\mathrm{pt}$.

The markings on the smooth locus of the stable core are parametrized as follows. Let $\mathscr F_\Gamma\to S_\Gamma$ be the family of smooth loci: pull back the universal curves of the factors $\mathcal M_{g_v,d_v}$ to $S_\Gamma$, remove the flag sections, and take the disjoint union over $v$. Let
\[
 X_{\Gamma,r}=(\mathscr F_\Gamma)^r_{S_\Gamma}\quad(r\ge1),\qquad
 X_{\Gamma,0}=S_\Gamma,
\]
the $r$th fiber power, which places $r$ labeled, possibly coincident, markings on the smooth core. The group $\Aut\Gamma$ acts on $S_\Gamma$ and $\mathscr F_\Gamma$ by permuting vertices and flags, and diagonally on $X_{\Gamma,r}$. For $\sigma\in\Aut\Gamma$, let
\begin{equation}\label{eq:vertex-series}
 \mathcal B_{\Gamma,\sigma}(w;u,v)=\sum_{r\ge0}
 E_\sigma(X_{\Gamma,r};u,v)\frac{w^r}{r!}.
\end{equation}
Only assignments of markings to vertices fixed by $\sigma$ contribute to the trace, since the remaining summands are permuted by $\sigma$. For example, if $\sigma$ reverses an edge joining two distinct vertices, neither vertex can carry a marking in a contributing assignment, whereas markings on a vertex with a reversed self-loop may contribute.

For an edge fixed by $\sigma$, the chain factor depends on whether $\sigma$ preserves or exchanges the two half-edges. Let $\iota$ be the involution of $\mathrm{LM}_k$ exchanging the endpoints and fixing every light label, and let
\[
 \mathcal C(w,uv)=\sum_{k\ge0}
 E(\mathrm{LM}_k;u,v)\frac{w^k}{k!},\qquad
 \mathcal D(w,uv)=\sum_{k\ge0}
 E_\iota(\mathrm{LM}_k;u,v)\frac{w^k}{k!}.
\]
Both depend on $(u,v)$ only through $uv$ (Lemma~\ref{lem:node-factors}).

\begin{proposition}\label{prop:graph-sum}
Let $g\ge2$. For a connected stable graph $\Gamma$ of genus $g$ without legs and $\sigma\in\Aut\Gamma$, let $p(\sigma)$ and $q(\sigma)$ be the numbers of edges fixed by $\sigma$ whose half-edges are, respectively, preserved and exchanged. Then
\begin{equation}\label{eq:graph-sum}
 \mathscr L_g(w;u,v)=
 \sum_{\Gamma}\frac1{|\Aut\Gamma|}
 \sum_{\sigma\in\Aut\Gamma}
 \mathcal B_{\Gamma,\sigma}(w;u,v)\,
 \mathcal C(w,uv)^{p(\sigma)}\,
 \mathcal D(w,uv)^{q(\sigma)},
\end{equation}
where $\Gamma$ runs over the finitely many isomorphism classes of such graphs.
\end{proposition}

\begin{proof}
Fix $\Gamma$ and $n$. The pullback of $\pi_n^{-1}(\mathcal M_\Gamma)$ to $S_\Gamma$ has a finite stratification by the assignments of labels to the smooth core and to the edges. For an assignment with $r$ labels on the core and $k_e$ labels at the edge $e$, the stratum is
\[
 X_{\Gamma,r}\times\prod_{e\in E(\Gamma)}\mathrm{LM}_{k_e},
\]
where, after ordering the half-edges of each edge, a nonempty chain is attached by identifying its ordered endpoints with the two flag points. No further gluing data are needed. The chains have no automorphisms preserving the endpoints and the labels, because each component has two attaching points and at least one marking, and automorphisms of the core components are already accounted for in $X_{\Gamma,r}$. The group $\Aut\Gamma$ permutes the strata, and the quotient of their union is $\pi_n^{-1}(\mathcal M_\Gamma)$.

Fix $\sigma\in\Aut\Gamma$ and compute its trace stratum by stratum. An edge moved by $\sigma$ carries no labels in a contributing stratum, so such edges contribute the factor $1$. On a fixed edge, $\sigma$ acts on the chain as the identity or as $\iota$, according to whether it preserves or exchanges the half-edges, giving the factors $\mathcal C$ and $\mathcal D$. The core contributes $\mathcal B_{\Gamma,\sigma}$. Since there are $n!/(r!\prod_ek_e!)$ assignments with given $r$ and $(k_e)$, multiplicativity and summation over $n$ show that $\sigma$ contributes $\mathcal B_{\Gamma,\sigma}\mathcal C^{p(\sigma)}\mathcal D^{q(\sigma)}$ to the exponential generating series. Averaging over $\Aut\Gamma$ by \eqref{eq:trace-average} and summing over the strata $\mathcal M_\Gamma$ of $\Mbar_g$ proves \eqref{eq:graph-sum}.
\end{proof}

\begin{example}[The theta graph]\label{ex:theta-core}
Put $t=uv$, and let $\Gamma_\theta$ consist of two genus-zero vertices $v_1,v_2$ joined by three edges. Then $S_{\Gamma_\theta}=\mathcal M_{0,3}^2$ is a point, the smooth locus of the core is $U=U_1\sqcup U_2$ with $U_i=\mathbb P^1\setminus\{0,1,\infty\}$, and $|\Aut\Gamma_\theta|=2\cdot3!=12$. Since $E(U)=2(t-2)$, the identity contributes $\mathcal B_{\Gamma_\theta,1}=e^{2(t-2)w}$, hence the term $e^{2(t-2)w}\mathcal C^3/12$.

If $\sigma$ fixes both vertices, then $E_\sigma(U_i)=1+t-\#\{\text{fixed flags at }v_i\}$, because the trace on the three deleted points counts fixed flags. The three edge transpositions therefore give $\mathcal B_{\Gamma_\theta,\sigma}=e^{2tw}$ with $(p,q)=(1,0)$, and the two $3$-cycles give $e^{2(t+1)w}$ with $(p,q)=(0,0)$. If $\sigma$ exchanges the vertices, it exchanges the two assignments of every marking, so only $r=0$ contributes and $\mathcal B_{\Gamma_\theta,\sigma}=1$. The exchange fixing every edge has $(p,q)=(0,3)$, the three exchanges fixing one edge have $(p,q)=(0,1)$, and the remaining two have $(p,q)=(0,0)$. Altogether, $\Gamma_\theta$ contributes
\[
 \frac1{12}\left(
 e^{2(t-2)w}\mathcal C^3
 +3e^{2tw}\mathcal C
 +2e^{2(t+1)w}
 +\mathcal D^3+3\mathcal D+2
 \right)
\]
to $\mathscr L_2(w;u,v)$, where $\mathcal C=\mathcal C(w,t)$ and $\mathcal D=\mathcal D(w,t)$.
\end{example}

We next compute the chain factors and show that the core factors are entire. Together, these two facts locate the poles of \eqref{eq:graph-sum}.

\begin{lemma}\label{lem:node-factors}
With $t=uv$,
\begin{equation}\label{eq:node-factors}
 \mathcal C(w,t)=\frac{t-1}{t-e^{(t-1)w}}=\frac1{1-W(w,t)},\qquad
 \mathcal D(w,t)=1+\frac{e^{(t+1)w}-1}{t+1}.
\end{equation}
In particular, $\mathcal D$ is entire in $(w,t)$.
\end{lemma}

\begin{proof}
The formula for $\mathcal C$ is the Poincar\'e series of Losev and Manin \cite[Theorem~2.3]{LosevManin}. Bergstr\"om and Minabe computed the $S_2\times S_k$-equivariant cohomology of $\mathrm{LM}_k$ \cite[Theorem~2.3]{BergstromMinabe}, and both formulas follow from \cite[Proposition~4.1]{BergstromMinabe} by specialization, as we now explain.

The cohomology of $\mathrm{LM}_k$ is concentrated in even degrees, with $H^{2j}$ of type $(j,j)$, so the variable $q$ of \cite{BergstromMinabe} is our $t$. In their notation, the power sums $p_m^x$ and $p_m^y$ record the cycle types of permutations of the endpoints and of the light labels. With the Frobenius characteristic normalized by $1/(2\cdot k!)$, the specialization $p_1^y=w$ and $p_m^y=0$ for $m\ge2$ selects the identity permutation of the labels and supplies the factor $w^k/k!$. In particular, the hook Schur functions specialize as
\[
 s_{(k-j,1^j)}^y\longmapsto \binom{k-1}{j}\frac{w^k}{k!},
\]
so the functions of \cite[Definition~2.2]{BergstromMinabe}, with $q=t$, satisfy
\[
 \begin{aligned}
 f_k^y&=\sum_{j=0}^{k-1}(-1)^j s_{(k-j,1^j)}^y t^{k-1-j}
 \longmapsto (t-1)^{k-1}\frac{w^k}{k!},\\
 g_k^y&=\sum_{j=0}^{k-1}s_{(k-j,1^j)}^y t^{k-1-j}
 \longmapsto (t+1)^{k-1}\frac{w^k}{k!}.
 \end{aligned}
\]
For the identity of $S_2$, set $(p_1^x)^2=2$ and $p_2^x=0$. The first term of \cite[(4.2)]{BergstromMinabe} becomes
\[
 \mathcal C(w,t)
 =\Bigl(1-\sum_{k\ge1}(t-1)^{k-1}\frac{w^k}{k!}\Bigr)^{-1}
 =\frac1{1-W(w,t)}.
\]
For the exchange of the endpoints, set $p_1^x=0$ and $p_2^x=2$. The second term of \cite[(4.2)]{BergstromMinabe} becomes
\[
 \Bigl(1+\sum_{k\ge1}g_k^y\Bigr)
 \Bigl(1-\sum_{k\ge1}p_2^y\circ f_k^y\Bigr)^{-1}.
\]
The plethysm $p_2^y\circ f_k^y$ replaces each $p_m^y$ by $p_{2m}^y$ and $t$ by $t^2$, so it vanishes under our specialization, and
\[
 \mathcal D(w,t)
 =1+\sum_{k\ge1}(t+1)^{k-1}\frac{w^k}{k!}
 =1+\frac{e^{(t+1)w}-1}{t+1}.
\]
The constant terms $1$ correspond to $\mathrm{LM}_0=\mathrm{pt}$. The singularity of $\mathcal D$ at $t=-1$ is removable, with value $1+w$.
\end{proof}

\begin{remark}\label{rem:chain-genus-zero}
Under the substitution \eqref{eq:substitution}, the chain factor $\mathcal C$ becomes the genus-zero series with two distinguished markings:
\[
  \mathcal C(T(z,t),t)
  =e^{T(z,t)}T_z(z,t)
  =\partial_z\varphi(z,t)
  =1+\sum_{n\ge1}E(\Mbar_{0,n+2};u,v)\frac{z^n}{n!},\qquad t=uv.
\]
Indeed, \eqref{eq:inverse-definition} and \eqref{eq:node-factors} give $R_w(w,t)=e^w/\mathcal C(w,t)$, so the first equality follows by differentiating $R(T(z,t),t)=z$. The others follow from \eqref{eq:genus-zero-log} and the definition of $\varphi$. Geometrically, the reduction morphism $\Mbar_{0,n+2}\to\mathrm{LM}_n$ 
contracts all components off the chain joining the two distinguished markings, which play the role of the two branches of a node.
\end{remark}

\begin{lemma}\label{lem:entire-vertices}
For $\sigma\in\Aut\Gamma$, the series $\mathcal B_{\Gamma,\sigma}$ in \eqref{eq:vertex-series} is entire in $w$, locally uniformly in $(u,v)$, and holomorphic in all variables.
\end{lemma}

\begin{proof}
Because of the factorial in \eqref{eq:vertex-series}, it suffices to bound $E_\sigma(X_{\Gamma,r};u,v)$ by $AB^r$, uniformly on compact sets of $(u,v)$, with $A,B$ independent of $r$. We first bound the Betti numbers.

Replace each positive-genus factor of $S_\Gamma$ by a finite \'etale cover with full level structure of level at least three, so that the product $S'$ is a smooth quasi-projective variety carrying universal curves. See \cite[Section~2.3.3]{PikaartDeJong} in genus at least two and \cite[Corollary~4.7.2]{KatzMazur} in genus one, taking one flag as the origin. Genus-zero factors need no cover, and the cover does not depend on $r$. Let $X'_{\Gamma,r}=X_{\Gamma,r}\times_{S_\Gamma}S'$. The universal smooth curves, together with their disjoint flag sections, are locally trivial as pairs in the analytic topology, so removing the sections and taking fiber powers produces locally trivial fiber bundles $X'_{\Gamma,r}\to S'$.

The variety $S'$ has the homotopy type of a finite CW complex. Let $K$ be its number of cells. For every local system $\mathscr A$ of finite rank, cellular cochains give
\begin{equation}\label{eq:cellular-local-system-bound}
 \sum_p\dim H^p(S';\mathscr A)\le K\rk\mathscr A,
\end{equation}
whatever the monodromy of $\mathscr A$. Let $B$ be the total Betti number of a fiber of $\mathscr F_\Gamma$, a disjoint union of punctured curves. By the K\"unneth formula, the higher direct images of the $r$th fiber power $f_r:X'_{\Gamma,r}\to S'$ have total rank $B^r$, so the Leray--Serre spectral sequence and \eqref{eq:cellular-local-system-bound} give
\begin{equation}\label{eq:exponential-betti-bound}
 \sum_i b_i(X'_{\Gamma,r})
 \le\sum_{p,q}\dim H^p(S';R^qf_{r*}\mathbb Q)
 \le KB^r.
\end{equation}
The cohomological transfer for the finite étale cover \(X'_{\Gamma,r}\to X_{\Gamma,r}\) shows that the cohomology of $X_{\Gamma,r}$ is a direct summand of that of $X'_{\Gamma,r}$, and Poincar\'e duality on these smooth spaces gives the same bound for compactly supported Betti numbers.

Finally, the eigenvalues of $\sigma$ on each graded piece have absolute value one, and the Hodge indices of $X_{\Gamma,r}$ are at most $\dim X_{\Gamma,r}=\dim S_\Gamma+r$. Hence, for every compact set $K_0\subset\mathbb C^2$, there are $A_{K_0},B_{K_0}$ with $\sup_{(u,v)\in K_0}|E_\sigma(X_{\Gamma,r};u,v)|\le A_{K_0}B_{K_0}^r$ for all $r$. This proves that the series in \eqref{eq:vertex-series} converges absolutely and uniformly on compact subsets. The level cover is used only for the Betti bound and need not be compatible with $\sigma$.
\end{proof}

\subsection{The pole of order \texorpdfstring{$3g-3$}{3g-3}}
Let $\mathscr G_g$ be the set of isomorphism classes of connected trivalent graphs of first Betti number $g$ without legs, with loops and multiple edges allowed and automorphisms acting on half-edges. Put
\begin{equation}\label{eq:core-constant}
 c_g=\sum_{\Gamma\in\mathscr G_g}\frac1{|\Aut\Gamma|}>0.
\end{equation}

\begin{proposition}\label{prop:hodge-pole}
For $g\ge2$ there is a function $\mathscr J_g(w;u,v)$, entire in $w$ locally uniformly in $(u,v)$ and holomorphic in all variables, such that, with $t=uv$,
\begin{equation}\label{eq:hodge-pole}
 \mathscr L_g(w;u,v)=
 \frac{c_g\, e^{2(g-1)(t-2)w}}{(1-W(w,t))^{3g-3}}
 +\frac{\mathscr J_g(w;u,v)}{(1-W(w,t))^{3g-4}}.
\end{equation}
In particular, for $t$ near $1$, $\mathscr L_g$ has a pole of order exactly $3g-3$ at $w=\tau(t)$.
\end{proposition}

\begin{proof}
By Lemmas~\ref{lem:node-factors} and~\ref{lem:entire-vertices}, each term of \eqref{eq:graph-sum} is an entire function divided by $(1-W)^{p(\sigma)}$, and $p(\sigma)\le|E(\Gamma)|=3g-3-\dim\mathcal M_\Gamma$ by \eqref{eq:stratum-dimension}. Equality $p(\sigma)=3g-3$ forces $\dim\mathcal M_\Gamma=0$, so that every vertex is rational and trivalent, and forces $\sigma$ to fix every half-edge, so that $\sigma=1$. For such $\Gamma$, the base $S_\Gamma$ is a point and the smooth locus of the corresponding core consists of $2g-2$ copies of $\mathbb P^1$ minus three points, with Hodge polynomial $2(g-1)(t-2)$. Hence $\mathcal B_{\Gamma,1}=e^{2(g-1)(t-2)w}$. Summing these terms gives the first term of \eqref{eq:hodge-pole}. All other terms have $p(\sigma)\le3g-4$ and can be put over the denominator $(1-W)^{3g-4}$ with an entire numerator. Since $W(\tau(t),t)=1$ and $W_w(\tau(t),t)=t\ne0$, the function $1-W$ has a simple zero at $\tau(t)$. The first numerator does not vanish there, and the second term has pole order at most $3g-4$, so the pole has order exactly $3g-3$.
\end{proof}

The leading term of \eqref{eq:hodge-pole} depends only on $uv$, so both specializations in \eqref{eq:parity-extraction} share it, and
\begin{equation}\label{eq:even-pole}
 \mathcal L_g(w,t)=
 \frac{c_g \,e^{2(g-1)(t-2)w}}{(1-W(w,t))^{3g-3}}
 +\frac{\mathcal{J}_g(w,t)}{(1-W(w,t))^{3g-4}},
\end{equation}
where $\mathcal{J}_g$ is entire in $w$, locally uniformly in $t$. We will use this for $t$ near $1$ and for $t$ on compact arcs of the unit circle.

\begin{remark}[Genus one as the case of exponent zero]\label{rem:genus-comparison}
In all genera, the leading singularity comes from rational curves whose dual graphs have only bivalent and trivalent vertices. Such a graph has $2g-2$ trivalent vertices. For $g\ge2$, contracting the bivalent chains produces a trivalent core with $3g-3$ edges, each contributing a factor $(1-W)^{-1}$. For $g=1$ there are no trivalent vertices, the graph is a cycle, and summing over cycles gives $\frac12\sum_{r\ge1}W^r/r=-\frac12\log(1-W)$, up to entire corrections from cycles with one or two components. The genus-one logarithm is thus the counterpart of the pole of order $3g-3$.
\end{remark}

\begin{proposition}\label{prop:core-constants}
The constants $c_g$ are determined by
\[
 \sum_{g\ge2}c_g x^{g-1}
 =\log\!\left(\sum_{m\ge0}
 \frac{(6m)!}{2^{3m}(3m)!\,(3!)^{2m}(2m)!}\,x^m\right).
\]
\end{proposition}

\begin{proof}
A trivalent graph with $2m$ vertices has $3m$ edges, hence genus $m+1$ if it is connected. Attach three half-edges to each of $2m$ labeled vertices. There are $(6m)!/(2^{3m}(3m)!)$ ways to pair the $6m$ half-edges, each pairing defining a graph, possibly with loops and multiple edges. The group $G_m$ permuting the vertices and, independently, the three half-edges at each vertex has order $(2m)!(3!)^{2m}$. Its orbits are the isomorphism classes of trivalent graphs with $2m$ vertices, and the stabilizer of a pairing representing $\Gamma$ is $\Aut\Gamma$. By the orbit--stabilizer theorem, dividing the number of pairings by $|G_m|$ gives $\sum_\Gamma1/|\Aut\Gamma|$ over all, possibly disconnected, such graphs. A graph consisting of $r$ copies of a connected graph $\Gamma$ has $r!|\Aut\Gamma|^r$ automorphisms, so the weighted series of all graphs is the exponential of that of connected graphs, the term $m=0$ corresponding to the empty graph. Taking logarithms proves the claim.
\end{proof}

\begin{example}[Low genus]\label{ex:core-constants}
Expanding the logarithm gives
\[
 c_2=\frac5{24},\qquad c_3=\frac5{16},\qquad
 c_4=\frac{1105}{1152},\qquad c_5=\frac{565}{128}.
\]
In genus two, $\mathscr G_2$ consists of the theta graph and the dumbbell graph:
\begin{center}
\begin{tikzpicture}[line width=0.6pt,vertex/.style={circle,fill,inner sep=0pt,minimum size=4pt},every label/.style={font=\small}]
 \begin{scope}
  \coordinate (a) at (0,0);
  \coordinate (b) at (2,0);
  \draw (a) -- (b);
  \draw (a) .. controls (0.5,1) and (1.5,1) .. (b);
  \draw (a) .. controls (0.5,-1) and (1.5,-1) .. (b);
  \node[vertex,label=above:$0$] at (a) {};
  \node[vertex,label=above:$0$] at (b) {};
  \node at (1,-1.5) {$|\Aut\Gamma|=12$};
 \end{scope}
 \begin{scope}[xshift=4.5cm]
  \coordinate (c) at (0,0);
  \coordinate (d) at (2,0);
  \draw (c) -- (d);
  \draw (c) .. controls (-1.3,1.6) and (-1.3,-1.6) .. (c);
  \draw (d) .. controls (3.3,1.6) and (3.3,-1.6) .. (d);
  \node[vertex,label=above:$0$] at (c) {};
  \node[vertex,label=above:$0$] at (d) {};
  \node at (1,-1.5) {$|\Aut\Gamma|=8$};
 \end{scope}
 \coordinate (corebox-sw) at (current bounding box.south west);
 \coordinate (corebox-ne) at (current bounding box.north east);
 \pgfresetboundingbox
 \path[use as bounding box] (corebox-sw) rectangle ([yshift=-7mm]corebox-ne);
\end{tikzpicture}
\end{center}
The labels indicate the vertex genera. The automorphism groups have orders $3!\cdot2=12$ and $2^2\cdot2=8$, so $c_2=\frac1{12}+\frac18=\frac5{24}$.
\end{example}
\section{Coefficient asymptotics and Gaussian limits}\label{sec:clt}

In this section, we extract uniform coefficient asymptotics from the singular expansions of Section~\ref{sec:light} and deduce the central limit theorems for the even and full Betti distributions, together with the asymptotics of the Euler characteristics.

\subsection{Asymptotics of the even Betti polynomials}\label{subsec:singularity-coefficients}

\begin{theorem}\label{thm:coefficients}
For every fixed $g\ge1$, as $n\to\infty$,
\begin{equation}\label{eq:coefficient-asymptotics}
 \begin{split}
 L_{g,n}(t)&=A_g(t)\,n!\,n^{3g-4}\,\tau(t)^{-n}
                      \bigl(1+O(n^{-1})\bigr),\\
 H_{g,n}(t)&=C_g(t)\,n!\,n^{\frac{3g-5}{2}}\rho(t)^{-n}
                      \bigl(1+O(n^{-\frac12})\bigr),
 \end{split}
\end{equation}
uniformly for complex $t$ in a fixed neighborhood of $1$. Here $A_1(t)=1/2$ and $C_1(t)=1/4$, and for $g\ge2$,
\begin{equation}\label{eq:amplitudes}
 A_g(t)=\frac{c_g\, e^{2(g-1)(t-2)\tau(t)}}{\Gamma(3g-3)\,[t\tau(t)]^{3g-3}},
 \qquad
 C_g(t)=\frac{c_g\, e^{\frac{g-1}{2}(t-2)\tau(t)}}{\Gamma\!\left(\frac{3g-3}{2}\right)[2\rho(t)]^{\frac{3g-3}{2}}}.
\end{equation}
The amplitudes are holomorphic and nonvanishing near $t=1$, and $A_g(1),C_g(1)>0$.
\end{theorem}

\begin{proof}
In each case we determine the singular expansion at the dominant singularity and apply the transfer theorem \cite[Theorem~VI.4]{FS} uniformly in $t$, as explained in Section~\ref{subsec:singularity-analysis}.

\emph{Light series.} To isolate the zero of $1-W(w,t)$ at $w=\tau(t)$, write
\[
 G(w,t)=\frac{1-W(w,t)}{1-w/\tau(t)}
 =t\tau(t)\,\frac{e^{(t-1)(w-\tau(t))}-1}{(t-1)(w-\tau(t))}.
\]
The singularity of $G$ at $w=\tau(t)$ is removable, with value $t\tau(t)$, and $G(w,1)=1$. Fix $R>1$. After shrinking the neighborhood of $t=1$, the function $G$ is holomorphic and bounded away from zero on a neighborhood of $|w|\le R$, its logarithm normalized at $w=0$ is uniformly bounded there, and $\sup_t|\tau(t)|<R$. We may therefore choose a fixed $\Delta$-domain $\Delta^L_0$ at $1$ with $\tau(t)\Delta^L_0\subset\{|w|<R\}$ for all such $t$. The light series continue analytically to $\tau(t)\Delta^L_0$.

For $g=1$, Lemma~\ref{lem:light-genus-one} gives the uniform expansion
\begin{equation}\label{eq:light-logarithmic-expansion}
 \mathcal L_1(w,t)=-\frac12\log\!\left(1-\frac{w}{\tau(t)}\right)
 +B(\tau(t),t)+O\!\left(1-\frac{w}{\tau(t)}\right),
\end{equation}
where $B(w,t)=\mathcal{J}_1(w,t)-\frac12\log G(w,t)$ is holomorphic and uniformly bounded on the disc. Transfer gives $[w^n]\mathcal L_1(w,t)=\tau(t)^{-n}\bigl(\frac1{2n}+O(n^{-2})\bigr)$, which is the first line of \eqref{eq:coefficient-asymptotics} with $A_1=1/2$. For $g\ge2$, the function $\mathcal L_g(w,t)$ in \eqref{eq:even-pole} is meromorphic on the disc with its only pole at $\tau(t)$, and
\begin{equation}\label{eq:light-pole-expansion}
 \mathcal L_g(w,t)=
 \frac{c_g\, e^{2(g-1)(t-2)\tau(t)}}{[t\tau(t)]^{3g-3}}
 \left(1-\frac{w}{\tau(t)}\right)^{-(3g-3)}
 +O\!\left(\left(1-\frac{w}{\tau(t)}\right)^{-(3g-4)}\right),
\end{equation}
so transfer gives the first line of \eqref{eq:coefficient-asymptotics} with relative error $O(n^{-1})$.

\emph{Heavy series.} We first check that composition with $T$ creates no new singularities on the domain $\Delta_t=\rho(t)\Delta_0$ of Lemma~\ref{lem:inverse}. We have
\begin{equation}\label{eq:matching}
 R_w(w,t)=e^w(1-W(w,t)),\qquad
 W(\tau(t),t)=1,\qquad
 R_w(T(z,t),t)\,T_z(z,t)=1.
\end{equation}
Since $T$ is holomorphic on $\Delta_t$, the last identity shows that $R_w(T(z,t),t)\ne0$, and hence $1-W(T(z,t),t)\ne0$, throughout $\Delta_t$. Next, since $W_w(\tau(t),t)=t$, Taylor expansion at $w=\tau(t)$ and Lemma~\ref{lem:inverse} give
\begin{equation}\label{eq:composed-denominator}
 \begin{split}
 1-W(T(z,t),t)&=-t\bigl(T(z,t)-\tau(t)\bigr)
 +O\!\left((T(z,t)-\tau(t))^2\right)\\
 &=(2\rho(t))^{\frac12}e^{\frac12(t-2)\tau(t)}
 \left(1-\frac{z}{\rho(t)}\right)^{\frac12}\left(1+O\!\left(\left(1-\frac{z}{\rho(t)}\right)^{\frac12}\right)\right).
 \end{split}
\end{equation}
For $g=1$, substituting \eqref{eq:composed-denominator} into \eqref{eq:genus-one-series} gives
\begin{equation}\label{eq:logarithmic-expansion}
 \mathcal H_1(z,t)=-\frac14\log\!\left(1-\frac{z}{\rho(t)}\right)
 +B_0(t)+O\!\left(\left(1-\frac{z}{\rho(t)}\right)^{\frac12}\right),
\end{equation}
with $B_0(t)=\mathcal{J}_1(\tau(t),t)-\frac14\log(2\rho(t))-\frac14(t-2)\tau(t)$ holomorphic near $t=1$. Transfer gives $[z^n]\mathcal H_1(z,t)=\rho(t)^{-n}\bigl(\frac1{4n}+O(n^{-3/2})\bigr)$, which is the second line of \eqref{eq:coefficient-asymptotics} with $C_1=1/4$. For $g\ge2$, substituting \eqref{eq:composed-denominator} into \eqref{eq:even-pole} via \eqref{eq:substitution} gives
\begin{equation}\label{eq:heavy-singular}
 \mathcal H_g(z,t)=
 \frac{c_g\, e^{\frac{g-1}{2}(t-2)\tau(t)}}{[2\rho(t)]^{\frac{3g-3}{2}}}
 \left(1-\frac{z}{\rho(t)}\right)^{-\frac{3g-3}{2}}
 +O\!\left(\left(1-\frac{z}{\rho(t)}\right)^{-\frac{3g-4}{2}}\right),
\end{equation}
and transfer gives the second line of \eqref{eq:coefficient-asymptotics} with relative error $O(n^{-1/2})$. The amplitudes \eqref{eq:amplitudes} are positive at $t=1$, hence nonvanishing nearby.
\end{proof}

\subsection{Central limit theorems}\label{subsec:central-euler}
Let $X^H_{g,n}$ and $X^L_{g,n}$ be the random variables with PGFs $H_{g,n}(t)/H_{g,n}(1)$ and $L_{g,n}(t)/L_{g,n}(1)$. They record half the degree of a random even-degree class. By Theorem~\ref{thm:coefficients}, their PGFs have the quasi-powers form, so Theorem~\ref{thm:quasi-powers} applies. The means are exactly $d_n/2$ by Poincar\'e duality.

\begin{definition}\label{def:variance-constant}
We set
\[
 \nu=\frac{3-e}{6(e-2)}=0.0653685\ldots,\qquad r_H=\frac1{4\nu}=\frac{3(e-2)}{2(3-e)}=3.824470\ldots.
\]
\end{definition}

The constant $\nu$ also governs the even Betti distributions of $\Mbar_{0,n}$ \cite[Theorem~3.1]{CKPNAS} and of the Fulton--MacPherson spaces of smooth projective curves \cite[Section~4]{CKTJM}. The constant $r_H$ will be the ultra-log-concavity threshold of $\Mbar_{g,n}$ (Theorem~\ref{thm:heavy-ulc}).

\begin{theorem}\label{thm:clt}
Fix $g\ge1$ and let $(\sigma^J_{g,n})^2=\Var(X^J_{g,n})$ for $J\in\{H,L\}$. Then $\mathbb E X^H_{g,n}=\mathbb E X^L_{g,n}=d_n/2$,
\[
 (\sigma^H_{g,n})^2=\nu n+\frac{(g-1)(11-4e)}{12(e-2)}+O(n^{-\frac12}),\qquad
 (\sigma^L_{g,n})^2=\frac{n-g+1}{12}+O(n^{-1}),
\]
and, for $J\in\{H,L\}$,
\[
 \sup_{x\in\mathbb R}\left|
 \mathbb P\left(\frac{X^J_{g,n}-d_n/2}{\sigma^J_{g,n}}\le x\right)
 -\Phi(x)\right|=O(n^{-\frac12})
\]
where $\Phi(x)$ denotes the standard normal distribution function.
\end{theorem}

\begin{proof}
Dividing the estimates of Theorem~\ref{thm:coefficients} by their values at $t=1$ gives
\begin{equation}\label{eq:normalized-pgfs}
 \begin{split}
 \frac{H_{g,n}(t)}{H_{g,n}(1)}
 &=\frac{C_g(t)}{C_g(1)}
   \left(\frac{e-2}{\rho(t)}\right)^n
   \bigl(1+O(n^{-\frac12})\bigr),\\
 \frac{L_{g,n}(t)}{L_{g,n}(1)}
 &=\frac{A_g(t)}{A_g(1)}
   \left(\frac1{\tau(t)}\right)^n
   \bigl(1+O(n^{-1})\bigr),
 \end{split}
\end{equation}
uniformly near $t=1$. The bases have the Taylor expansions
\[
 \begin{aligned}
 \frac{e-2}{\rho(t)}&=1+\frac{t-1}{2}-\frac{5e-12}{24(e-2)}(t-1)^2
                  +O\bigl((t-1)^3\bigr),\\
 \frac1{\tau(t)}&=1+\frac{t-1}{2}-\frac{(t-1)^2}{12}
                  +O\bigl((t-1)^3\bigr),
 \end{aligned}
\]
so that
\[
 \mathfrak m\left(\frac{e-2}{\rho}\right)
 =\mathfrak m\left(\frac1{\tau}\right)=\frac12,\qquad
 \mathfrak v\left(\frac{e-2}{\rho}\right)=\nu>0,\qquad
 \mathfrak v\left(\frac1{\tau}\right)=\frac1{12}>0.
\]
For $g\ge2$, taking logarithms in \eqref{eq:amplitudes} and using $(t-2)\tau(t)=\log t-\tau(t)$ gives
\begin{equation}\label{eq:amplitude-logarithms}
 \begin{split}
 \log C_g(t)&=\mathrm{const}
       +\frac{g-1}{2}\bigl(\log t-\tau(t)-3\log\rho(t)\bigr),\\
 \log A_g(t)&=\mathrm{const}
       -(g-1)\bigl(\log t+2\tau(t)+3\log\tau(t)\bigr),
 \end{split}
\end{equation}
and differentiating at $t=1$ yields
\[
 \mathfrak m(C_g)=\mathfrak m(A_g)=\frac{3g-3}{2},\qquad
 \mathfrak v(C_g)=\frac{(g-1)(11-4e)}{12(e-2)},
 \qquad \mathfrak v(A_g)=-\frac{g-1}{12}.
\]
For $g=1$ the amplitudes are constant and these quantities vanish. Theorem~\ref{thm:quasi-powers}, applied with $\beta_n=n$ and $\varepsilon_n=n^{-1/2}$ in the heavy case and $\varepsilon_n=n^{-1}$ in the light case, completes the proof.
\end{proof}

\begin{corollary}\label{cor:full-degree}
For each fixed $g\ge1$, the full Betti distributions of $\Mbar_{g,n}$ and $\Mbar_{g,0|n}$, indexed by cohomological degree, are asymptotically normal. They have exact mean $d_n=n+3g-3$, their variances are given by four times the expansions in Theorem~\ref{thm:clt}, with the same error orders. After centering and scaling their distribution functions converge uniformly to $\Phi$ at rate $O(n^{-1/2})$.
\end{corollary}

\begin{proof}
By \eqref{eq:parity-extraction}, the Poincar\'e polynomials are obtained by setting $u=v=-q$ in the Hodge--Deligne polynomials. The leading terms in \eqref{eq:genus-one-hodge-series} and \eqref{eq:hodge-pole} depend only on $uv=q^2$, so the odd-degree contributions are contained in the remainders. Using \eqref{eq:substitution} and the argument of Theorem~\ref{thm:coefficients}, the Poincar\'e polynomials therefore have the same leading terms and relative errors as the even Betti polynomials at $t=q^2$, uniformly for $q$ near $1$. Normalizing at $q=1$ gives the quasi-powers forms \eqref{eq:normalized-pgfs} with $t=q^2$. Since $\mathfrak m(f(q^2))=2\mathfrak m(f)$ and $\mathfrak v(f(q^2))=4\mathfrak v(f)$ for $f$ holomorphic near $1$ with $f(1)\ne0$, Theorem~\ref{thm:quasi-powers} gives the claim. The means are given by Poincar\'e duality.
\end{proof}

\begin{figure}[t]
\centering
\begin{minipage}{0.48\textwidth}
\centering
\includegraphics[width=\linewidth]{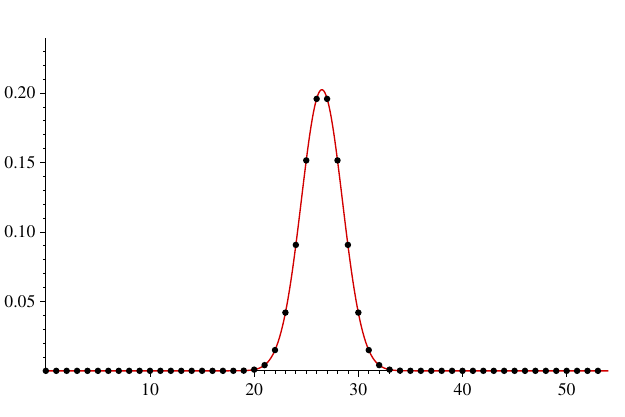}
\par\smallskip
$\Mbar_{2,50}$
\end{minipage}\hfill
\begin{minipage}{0.48\textwidth}
\centering
\includegraphics[width=\linewidth]{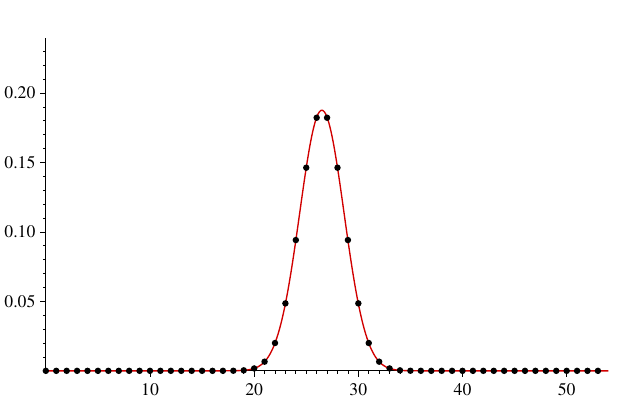}
\par\smallskip
$\Mbar_{2,0|50}$
\end{minipage}
\caption{Even Betti distributions for $g=2$ and $n=50$. The horizontal coordinate is the half-degree $k$. Black dots show the exact normalized even Betti numbers $b_{2k}/\sum_j b_{2j}$, and red curves show the Gaussian densities with the exact means and variances.}
\label{fig:genus-two-even}
\end{figure}

\begin{example}\label{ex:genus-two-figure}
Figure~\ref{fig:genus-two-even} compares the exact normalized even Betti numbers of $\Mbar_{2,50}$ and $\Mbar_{2,0|50}$ with the Gaussian densities having the same means and variances. Both distributions have mean $53/2$, and their variances are approximately $3.885$ and $4.521$. At the middle indices $k=26,27$, the Gaussian approximations have relative errors of approximately $0.0913\%$ and $0.1426\%$.
\end{example}

\subsection{Euler characteristics}\label{subsec:euler-ratio}
Kannan, Serpente and Yun asked for the asymptotics of $\chi(\Mbar_{g,0|n})/\chi(\Mbar_{g,n})$ for fixed $g\ge2$, and in particular whether this ratio tends to zero \cite[Question~5.3]{KSY}. For $g=1$, they proved that it is asymptotic to $2(e-2)^n$ \cite[Corollary~5.2]{KSY}. Here $\chi$ is the Euler characteristic of the coarse space, equivalently the alternating sum of the rational Betti numbers of the stack.

\begin{corollary}\label{cor:euler-ratio}
For every fixed $g\ge2$,
\begin{equation}\label{eq:euler-ratio}
 \frac{\chi(\Mbar_{g,0|n})}{\chi(\Mbar_{g,n})}
 =\left(\frac{2(e-2)}{e}\right)^{\frac{3g-3}{2}}
 \frac{\Gamma\!\left(\frac{3g-3}{2}\right)}{\Gamma(3g-3)}\,
 n^{\frac{3g-3}{2}}(e-2)^n\bigl(1+O(n^{-\frac12})\bigr).
\end{equation}
In particular, the ratio tends to zero. For $g=1$, the same argument recovers the asymptotic $2(e-2)^n$.
\end{corollary}

\begin{proof}
The estimates in the proof of Corollary~\ref{cor:full-degree}, evaluated at $q=1$, compare the sum of all Betti numbers with the sum of the even Betti numbers.
\begin{equation}\label{eq:full-even-ratios}
 \frac{\sum_{j}b_j(\Mbar_{g,n})}{H_{g,n}(1)}=1+O(n^{-\frac12}), \qquad
 \frac{\sum_{j}b_j(\Mbar_{g,0|n})}{L_{g,n}(1)}=1+O(n^{-1}).
\end{equation}
Since the Euler characteristic is twice the even Betti total minus the total Betti number, \begin{equation}
    \chi(\Mbar_{g,n})=H_{g,n}(1)(1+O(n^{-1/2})) \quad \text{and}\quad \chi(\Mbar_{g,0|n})=L_{g,n}(1)(1+O(n^{-1})).\end{equation}
    Theorem~\ref{thm:coefficients} at $t=1$ now gives
\[
 \frac{\chi(\Mbar_{g,0|n})}{\chi(\Mbar_{g,n})}
 =\frac{A_g(1)}{C_g(1)}\,n^{\frac{3g-3}{2}}(e-2)^n
       \bigl(1+O(n^{-\frac12})\bigr),
 \qquad
 \frac{A_g(1)}{C_g(1)}
   =\left(\frac{2(e-2)}{e}\right)^{\frac{3g-3}{2}}
   \frac{\Gamma\!\left(\frac{3g-3}{2}\right)}{\Gamma(3g-3)},
\]
which is \eqref{eq:euler-ratio}. The ratio tends to zero because $0<e-2<1$. For $g=1$, the same computation with $A_1/C_1=2$ gives $2(e-2)^n(1+O(n^{-1/2}))$.
\end{proof}

By \eqref{eq:full-even-ratios}, odd cohomology forms a proportion $O(n^{-1/2})$ of the cohomology of $\Mbar_{g,n}$ and $O(n^{-1})$ of that of $\Mbar_{g,0|n}$. Theorem~\ref{thm:hodge-asymmetry} below shows that these orders are sharp for $g\ge2$ and determines the leading constants.

\section{Local limits and central ultra-log-concavity}\label{sec:local}

The central limit theorems of Section~\ref{sec:clt} control cumulative probabilities but not individual Betti numbers. In this section, we prove a local limit theorem with Edgeworth-type corrections, deduce asymptotic formulas for the individual even Betti numbers, and then compare neighboring Betti numbers precisely enough to determine the central ultra-log-concavity thresholds. The local limit theorem is a refinement of Hwang's Fourier method \cite[Theorem~1 and Section~2]{HwangLLT} that exploits the symmetry provided by Poincar\'e duality.

\subsection{A symmetric local limit theorem}\label{subsec:symmetric-local}
A Gaussian approximation to individual probabilities suffices for the Betti number asymptotics, but ultra-log-concavity compares neighboring probabilities and requires more precision. Since the distributions are symmetric, the odd centered cumulants vanish. We therefore establish a symmetric version of the local limit theorem with higher-order corrections involving only even cumulants. As usual, $\He_j(z)=(-1)^je^{z^2/2}\frac{d^j}{dz^j}e^{-z^2/2}$ denotes the probabilists' Hermite polynomial of degree $j$.

\begin{theorem}[Symmetric quasi-powers local limit theorem]
\label{thm:symmetric-local}
Let $(Y_n)$ be integer-valued random variables with moment generating functions $M_n(s)=\mathbb E e^{sY_n}$, and assume the following.
\begin{enumerate}
\item \emph{Quasi-powers.} There are $r_0>0$, a sequence $\varepsilon_n\to0$, and functions $\lambda,\psi$ independent of $n$, holomorphic on a neighborhood of $|s|\le r_0$ and real on the real axis, with $\lambda(0)=\psi(0)=0$ and $\lambda''(0)>0$, such that $M_n$ is holomorphic there and
\begin{equation}\label{eq:symmetric-quasi}
 M_n(s)=\exp\bigl(n\lambda(s)+\psi(s)\bigr)
               \bigl(1+O(\varepsilon_n)\bigr)
\end{equation}
uniformly for complex $|s|\le r_0$.
\item \emph{Symmetry.} With $\mu_n=\mathbb E Y_n$, we have $\mathbb P(Y_n=\mu_n+x)=\mathbb P(Y_n=\mu_n-x)$ for all $x\in\mathbb R$.
\item \emph{Decay away from zero.} For every sufficiently small $\delta>0$, there is $c_\delta>0$ with
\begin{equation}\label{eq:symmetric-off-zero}
 \sup_{\delta\le|\theta|\le\pi}|M_n(i\theta)|
       =O(e^{-c_\delta n}).
\end{equation}
\end{enumerate}
Write $p_{n,k}=\mathbb P(Y_n=k)$, $\sigma_n^2=\Var(Y_n)$, $\kappa_{j,n}=\kappa_j(Y_n)$ and $z_{n,k}=(k-\mu_n)/\sigma_n$. Then
\begin{equation}\label{eq:symmetric-cumulants}
 \sigma_n^2=n\lambda''(0)+\psi''(0)+O(\varepsilon_n),\qquad
 \kappa_{j,n}=n\lambda^{(j)}(0)+\psi^{(j)}(0)+O(\varepsilon_n)
                                      \quad (j=4,6),
\end{equation}
and, uniformly for all $k\in\mathbb Z$,
\begin{equation}\label{eq:local-sixth}
 \sigma_np_{n,k}=
 \frac{e^{-z_{n,k}^2/2}}{\sqrt{2\pi}}\bigg[
  1+\frac{\kappa_{4,n}}{4!\,\sigma_n^4}\He_4(z_{n,k})
  +\frac{\kappa_{6,n}}{6!\,\sigma_n^6}\He_6(z_{n,k}) +\frac{\kappa_{4,n}^2}{2(4!)^2\sigma_n^8}\He_8(z_{n,k})
 \bigg]+O(n^{-3}).
\end{equation}
\end{theorem}

\begin{proof}
The cumulant estimates \eqref{eq:symmetric-cumulants} are \eqref{eq:cumulants}. For the local expansion, fix a small $\delta\in(0,r_0/2)$ and split the Fourier integral \eqref{eq:prelim-fourier} as $p_{n,k}=I_1+I_2$, with
\[
 I_1=\frac1{2\pi}\int_{-\delta}^{\delta}
                       e^{-ik\theta}M_n(i\theta)\,d\theta,
 \qquad
 I_2=\frac1{2\pi}\int_{\delta<|\theta|\le\pi}
                       e^{-ik\theta}M_n(i\theta)\,d\theta.
\]
By condition~(3), $|I_2|=O(e^{-c_\delta n})$ uniformly in $k$, so the expansion comes from $I_1$.

Let $K_n(s)=\log M_n(s)-\mu_ns$ be the centered cumulant generating function. By \eqref{eq:cumulants} with $j=8$, we have $K_n^{(8)}(s)=O(n)$ for $|s|\le r_0/2$, and $K_n$ is even by condition~(2). Taylor's formula therefore gives
\begin{equation}\label{eq:probability-log-taylor}
 K_n(s)=\frac{\sigma_n^2s^2}{2}
        +\frac{\kappa_{4,n}s^4}{4!}
        +\frac{\kappa_{6,n}s^6}{6!}+O(n|s|^8),\qquad |s|\le r_0/2.
\end{equation}
Substituting $s=i\tau/\sigma_n$, we obtain
\[
 \sigma_nI_1=\frac1{2\pi}\int_{-\delta\sigma_n}^{\delta\sigma_n}
             e^{-i\tau z_{n,k}}e^{K_n(i\tau/\sigma_n)}\,d\tau,
 \qquad
 K_n(i\tau/\sigma_n)=-\frac{\tau^2}2
  +\frac{\kappa_{4,n}\tau^4}{4!\,\sigma_n^4}
  -\frac{\kappa_{6,n}\tau^6}{6!\,\sigma_n^6}
 +O(n^{-3}\tau^8).
\]
Fix $0<\alpha<1/4$ and consider first the range $|\tau|\le n^\alpha$. There the terms after $-\tau^2/2$ are $o(1)$ uniformly, since $\tau^4/n\le n^{4\alpha-1}$, and expanding their exponential gives
\begin{equation}\label{eq:probability-fourier-error}
 e^{K_n(i\tau/\sigma_n)}
 =e^{-\tau^2/2}P_n(\tau)
 +O\!\left(n^{-3}e^{-\tau^2/2}
             (\tau^8+\tau^{10}+\tau^{12})\right),
\end{equation}
where
\[
  P_n(\tau)=1+\frac{\kappa_{4,n}\tau^4}{4!\,\sigma_n^4}
  -\frac{\kappa_{6,n}\tau^6}{6!\,\sigma_n^6}
  +\frac{\kappa_{4,n}^2\tau^8}{2(4!)^2\sigma_n^8}.
\]
The coefficients of $\tau^4$, $\tau^6$ and $\tau^8$ in $P_n$ are $O(n^{-1})$, $O(n^{-2})$ and $O(n^{-2})$, and the error term in \eqref{eq:probability-fourier-error} integrates to $O(n^{-3})$, uniformly in $k$ because $|e^{-i\tau z_{n,k}}|=1$.

On the range $n^\alpha<|\tau|\le\delta\sigma_n$, \eqref{eq:probability-log-taylor} and \eqref{eq:symmetric-cumulants} give $\Re K_n(i\tau/\sigma_n)=-\tau^2/2+O(\tau^4/n)$, and $\tau^4/n\le\delta^2(\sigma_n^2/n)\tau^2$. After reducing $\delta$, we therefore have $|e^{K_n(i\tau/\sigma_n)}|\le e^{-\tau^2/4}$ for $|\tau|\le\delta\sigma_n$, so this range contributes $O(e^{-cn^{2\alpha}})$ for some $c>0$. Extending the integral of $e^{-\tau^2/2}P_n(\tau)$ from $|\tau|\le n^\alpha$ to $\mathbb R$ costs $O(e^{-n^{2\alpha}/4})$, since the coefficients of $P_n$ are bounded. Hence
\[
 \sigma_nI_1=\frac1{2\pi}\int_{\mathbb R}
 e^{-i\tau z_{n,k}}e^{-\tau^2/2}P_n(\tau)\,d\tau
 +O(n^{-3}).
\]
Differentiating the Fourier transform of the Gaussian gives
\begin{equation}\label{eq:probability-gaussian-fourier}
 \frac1{2\pi}\int_{\mathbb R}
 e^{-i\tau z}\tau^{2j}e^{-\tau^2/2}\,d\tau
 =(-1)^j\He_{2j}(z)\frac{e^{-z^2/2}}{\sqrt{2\pi}}
 \qquad (j\ge0),
\end{equation}
so $\tau^4,\tau^6,\tau^8$ produce $\He_4,-\He_6,\He_8$ respectively. The sign of the sixth-cumulant term in $P_n$ cancels the sign of $\He_6$. Finally, $\sigma_nI_2=O(\sqrt n\,e^{-c_\delta n})=o(n^{-3})$, which proves \eqref{eq:local-sixth}.
\end{proof}

\subsection{Local limits for the even Betti numbers}\label{subsec:betti-local}
We verify the hypotheses of Theorem~\ref{thm:symmetric-local} for $X^H_{g,n}$ and $X^L_{g,n}$. Substituting $t=e^s$ in \eqref{eq:normalized-pgfs} gives
\[
 \mathbb E e^{sX^H_{g,n}}
 =e^{n\lambda_H(s)+\psi_H(s)}
   \bigl(1+O(n^{-\frac12})\bigr),\qquad
 \mathbb E e^{sX^L_{g,n}}
 =e^{n\lambda_L(s)+\psi_L(s)}
   \bigl(1+O(n^{-1})\bigr),
\]
uniformly on a fixed complex disc, where
\begin{equation}\label{eq:lambda-functions}
 \begin{aligned}
 \lambda_H(s)&=\log\frac{e-2}{\rho(e^s)},&\qquad
 \lambda_L(s)&=-\log\tau(e^s)
              =\log\frac{e^s-1}{s},\\
 \psi_H(s)&=\log\frac{C_g(e^s)}{C_g(1)},&\qquad
 \psi_L(s)&=\log\frac{A_g(e^s)}{A_g(1)},
 \end{aligned}
\end{equation}
all logarithms vanishing at $s=0$. By \eqref{eq:analytic-moments-chain} and the proof of Theorem~\ref{thm:clt},
\begin{equation}\label{eq:lambda-psi-derivatives}
 \begin{gathered}
 \lambda_H'(0)=\lambda_L'(0)=\frac12,\qquad
 \psi_H'(0)=\psi_L'(0)=\frac{3g-3}{2},\qquad
 \lambda_H''(0)=\nu,\qquad \lambda_L''(0)=\frac1{12},\\
 \psi_H''(0)=\frac{(g-1)(11-4e)}{12(e-2)},\qquad \psi_L''(0)=-\frac{g-1}{12}.
 \end{gathered}
\end{equation}
This is condition~(1), with $\varepsilon_n=n^{-1/2}$ and $\varepsilon_n=n^{-1}$, respectively. Poincar\'e duality gives the symmetry about $d_n/2$ required in condition~(2). Condition~(3) is the following proposition.

\begin{proposition}\label{prop:off-zero}
For every fixed $g\ge1$ and $\delta\in(0,\pi)$, there are $c_1,c_2>0$ such that
\begin{equation}\label{eq:off-zero}
 \sup_{\delta\le|\theta|\le\pi}
 \bigl|\mathbb E e^{i\theta X^H_{g,n}}\bigr|
 =O(e^{-c_1n}),\qquad
 \sup_{\delta\le|\theta|\le\pi}
 \bigl|\mathbb E e^{i\theta X^L_{g,n}}\bigr|
 =O(e^{-c_2n}).
\end{equation}
\end{proposition}

\begin{proof}
Write $t=e^{i\theta}$ with $\delta\le|\theta|\le\pi$ and $\log t=i\theta$.

\emph{Light case.} The zeros of $1-W(w,t)$ closest to the origin have modulus $|\tau(t)|=|\theta|/(2|\sin(\theta/2)|)>1$. Choose $1<R_L<\min_{\delta\le|\theta|\le\pi}|\tau(e^{i\theta})|$. By \eqref{eq:genus-one-series} and \eqref{eq:even-pole}, $\mathcal L_g(w,t)$ is holomorphic on a neighborhood of $|w|\le R_L$ and uniformly bounded on $|w|=R_L$ for $\delta\le|\theta|\le\pi$.

\emph{Heavy case.} Expanding \eqref{eq:inverse-definition} gives
\[
 R(w,t)=w-\sum_{m\ge3}t(1+t+\cdots+t^{m-3})\frac{w^m}{m!}.
\]
The coefficient polynomial has absolute value at most $m-2$, and for $m=4$ it satisfies the sharper bound $|t+t^2|=2|\cos(\theta/2)|\le2\cos(\delta/2)$. Since $\sum_{m\ge3}(m-2)/m!=3-e$, we get, for $|w|\le1$,
\[
 |R(w,t)-w|\le3-e-\frac{1-\cos(\delta/2)}{12}.
\]
Let $R_H=e-2+(1-\cos(\delta/2))/24>e-2$. For $|z|\le R_H$ and $|w|=1$, we have $|(R(w,t)-z)-w|\le1-(1-\cos(\delta/2))/24<1$. By Rouch\'e's theorem, the equation $R(w,t)=z$ has exactly one root in $|w|<1$, which is simple. By the implicit function theorem and uniqueness, these roots define a holomorphic continuation of $T(z,t)$ to a neighborhood of $|z|\le R_H$. Moreover,
\[
 |T(z,t)|\le |z|+|R(T(z,t),t)-T(z,t)|
 \le1-\frac{1-\cos(\delta/2)}{24},
\]
so $T(z,t)$ stays in a compact subset of the unit disc, where $\mathcal L_g(w,t)$ is holomorphic. By \eqref{eq:substitution}, $\mathcal H_g(z,t)$ is holomorphic on a neighborhood of $|z|\le R_H$ and uniformly bounded on $|z|=R_H$.

\emph{Conclusion.} Cauchy's estimates give $|H_{g,n}(e^{i\theta})|/n!=O(R_H^{-n})$ and $|L_{g,n}(e^{i\theta})|/n!=O(R_L^{-n})$ uniformly in $\theta$. By Theorem~\ref{thm:coefficients}, $H_{g,n}(1)/n!$ and $L_{g,n}(1)/n!$ are bounded below by positive multiples of $n^{-1}(e-2)^{-n}$ and $n^{-1}$. Hence the two suprema in \eqref{eq:off-zero} are $O(n((e-2)/R_H)^n)$ and $O(nR_L^{-n})$, which proves the proposition.
\end{proof}

\begin{theorem}\label{thm:llt}
Fix $g\ge1$ and $J\in\{H,L\}$. Put $p_{n,k}=\mathbb P(X^J_{g,n}=k)$, $\sigma_n^2=\Var(X^J_{g,n})$, $z_{n,k}=(k-d_n/2)/\sigma_n$, and let $\kappa_{j,n}$ be the $j$th cumulant of $X^J_{g,n}$. Then \eqref{eq:local-sixth} holds uniformly for all $k\in\mathbb Z$.
\end{theorem}

\begin{proof}
The three hypotheses of Theorem~\ref{thm:symmetric-local} were verified above.
\end{proof}

Multiplying the local probabilities by the even Betti totals of Theorem~\ref{thm:coefficients} gives the individual Betti numbers.

\begin{corollary}\label{cor:even-betti-estimates}
Fix $g\ge1$ and put $x_{n,k}=k-d_n/2$. Uniformly over all integers $k$,
\begin{equation}\label{eq:even-betti-estimates}
 \begin{split}
 b_{2k}(\Mbar_{g,n})
 &=\frac{C_g(1)}{\sqrt{2\pi\nu}}\,
 n!\,(e-2)^{-n}n^{\frac{3g-6}{2}}
 \left[\exp\!\left(-\frac{x_{n,k}^2}{2\nu n}\right)
 +O(n^{-\frac12})\right],\\
 b_{2k}(\Mbar_{g,0|n})
 &=\sqrt{\frac6\pi}\,A_g(1)\,
 n!\,n^{\frac{6g-9}{2}}
 \left[\exp\!\left(-\frac{6x_{n,k}^2}{n}\right)
 +O(n^{-1})\right],
 \end{split}
\end{equation}
where $A_1(1)=\frac12$, $C_1(1)=\frac14$, and for $g\ge2$
\[
 A_g(1)=\frac{c_g\, e^{-2(g-1)}}{\Gamma(3g-3)},\qquad
 C_g(1)=\frac{c_g\, e^{-\frac{g-1}{2}}}
 {\Gamma\!\left(\frac{3g-3}{2}\right)[2(e-2)]^{\frac{3g-3}{2}}}.
\]
\end{corollary}

\begin{proof}
Since $\kappa_{4,n},\kappa_{6,n}=O(n)$ and each $e^{-y^2/2}\He_j(y)$ is bounded on $\mathbb R$, Theorem~\ref{thm:llt} gives
\[
 \mathbb P(X^J_{g,n}=k)
 =\frac{1}{\sqrt{2\pi}\,\sigma^J_{g,n}}
   \left[\exp\!\left(-\frac{x_{n,k}^2}{2(\sigma^J_{g,n})^2}\right)
   +O(n^{-1})\right]
\]
uniformly in $k$. By Theorem~\ref{thm:clt}, $(\sigma^J_{g,n})^2=v_Jn+O(1)$ with $v_H=\nu$ and $v_L=1/12$. By the mean value theorem, replacing $(\sigma^J_{g,n})^2$ by $v_Jn$ in the exponent changes the Gaussian factor by $O(n^{-1})$ uniformly in $k$, and $1/\sigma^J_{g,n}=(v_Jn)^{-1/2}(1+O(n^{-1}))$. Multiplying by $H_{g,n}(1)$ or $L_{g,n}(1)$ from Theorem~\ref{thm:coefficients} and evaluating \eqref{eq:amplitudes} at $t=1$ completes the proof.
\end{proof}

\subsection{Central ultra-log-concavity}\label{subsec:central-ulc}
For a function $f$ on the integers, write $\Delta^2 f(k)=f(k-1)-2f(k)+f(k+1)$. For a positive sequence $(a_k)_{k=0}^{d}$ and $1\le k\le d-1$, $r$-ultra-log-concavity at $k$ means
\begin{equation}\label{eq:ulc-log-criterion}
 \Delta^2\log\left(\frac{a_k}{\binom{d}{k}^r}\right)\le0,
\end{equation}
and strict $r$-ultra-log-concavity means strict inequality.

The strategy is to compare $\Delta^2\log p_{n,k}$ with $r\Delta^2\log\binom{d_n}k$. Both are of order $n^{-1}$, and when $d_n=n+O(1)$ the leading terms agree exactly when $r$ equals $1/(4v)$, where $vn$ is the leading term of the variance. At this critical exponent the terms of order $n^{-1}$ cancel, and we must expand to order $n^{-2}$. The sixth-cumulant and squared fourth-cumulant terms in \eqref{eq:local-sixth} then contribute only $o(n^{-2})$ on the central range. Hence only the variance and the fourth cumulant matter contribute at order $n^{-2}$.

\begin{lemma}[Central ultra-log-concavity criterion]\label{lem:central-ulc}
Suppose that $(Y_n)$ satisfies the hypotheses of Theorem~\ref{thm:symmetric-local}, that $Y_n$ takes values in $\{0,\ldots,d_n\}$ with $d_n=n+O(1)$, and that $\mu_n=d_n/2$. Put $v=\lambda''(0)>0$, so that
\begin{equation}\label{eq:ulc-parameters}
 \sigma_n^2=vn+O(1).
\end{equation}
Fix $0<c_0<\sqrt2$. The following statements hold uniformly for all $k$ such that $|z_{n,k+j}|\le c_0\sqrt{\log n}$ for $j=-1,0,1$.
\begin{enumerate}
\item Put $r_*=1/(4v)$. For each fixed $0\le r<r_*$, the sequence $(p_{n,k})_k$ is strictly $r$-ultra-log-concave at $k$ for all sufficiently large $n$. For each fixed $r>r_*$, the strict reverse of \eqref{eq:ulc-log-criterion} holds at $k$ for all sufficiently large $n$.
\item The sequence $(p_{n,k})_k$ is strictly $r_*$-ultra-log-concave at $k$ for all sufficiently large $n$ if
\begin{equation}\label{eq:ulc-critical-condition}
 \lambda^{(4)}(0)\le-8v^3\qquad\text{and}\qquad \psi''(0)<v\bigl(2\psi'(0)+1\bigr)+\frac{\lambda^{(4)}(0)}{2v}.
\end{equation}
\end{enumerate}
\end{lemma}

\begin{proof}
Put $L_n=c_0\sqrt{\log n}$. On the range $|z_{n,k}|\le L_n$, the Gaussian factor $e^{-z_{n,k}^2/2}/\sqrt{2\pi}$ is at least $n^{-c_0^2/2}/\sqrt{2\pi}$, and $c_0^2/2<1$. A unit step in $k$ changes $z_{n,k}$ by $\sigma_n^{-1}=O(n^{-1/2})$.

\textup{(1)} Keeping only the fourth-cumulant term in \eqref{eq:local-sixth} gives
\[
 \sigma_np_{n,k}=\frac{e^{-z_{n,k}^2/2}}{\sqrt{2\pi}}
 \left[1+\frac{\kappa_{4,n}}{4!\,\sigma_n^4}\He_4(z_{n,k})\right]+O(n^{-2}).
\]
After division by the Gaussian factor, the remainder becomes $O(n^{-2}e^{L_n^2/2})$, while the correction term is $O(n^{-1}(1+L_n^4))$. Both tend to zero uniformly, so $p_{n,k}>0$ for large $n$, and
\[
 \log p_{n,k}=-\log(\sqrt{2\pi}\,\sigma_n)-\frac{z_{n,k}^2}2
        +\frac{\kappa_{4,n}}{4!\,\sigma_n^4}\He_4(z_{n,k})
        +O\!\left(n^{-2}(1+L_n^8+e^{L_n^2/2})\right).
\]
The second difference of the fourth-cumulant term is $O(n^{-2}(1+L_n^2))$, and that of the remainder is $O(n^{-2}e^{L_n^2/2})=o(n^{-1})$. Hence, by \eqref{eq:ulc-parameters},
\[
 \Delta^2\log p_{n,k}=-\sigma_n^{-2}+o(n^{-1})=-\frac1{vn}+o(n^{-1}).
\]
On the other hand, since $d_n=n+O(1)$ and $|k-d_n/2|\le\sigma_nL_n$,
\[
 \Delta^2\log\binom{d_n}k
 =\log\frac{k(d_n-k)}{(k+1)(d_n-k+1)}
 =-\frac4n+O\!\left(n^{-2}(1+L_n^2)\right).
\]
Therefore
\begin{equation}\label{eq:ulc-curvature}
 \Delta^2\log\left(\frac{p_{n,k}}{\binom{d_n}k^r}\right)
 =\frac{4r-v^{-1}}n+o_{r}(n^{-1}),
\end{equation}
and for $r\ne r_*$ the sign is that of $4r-v^{-1}$.

\textup{(2)} We now use the full expansion \eqref{eq:local-sixth}. The additional terms in $\log p_{n,k}$ involve $\He_6$, $\He_8$ and $\He_4^2$, with coefficients $O(n^{-2})$. As these polynomials have degree at most $8$, their second differences are $O(n^{-3}(1+L_n^8))$. The Taylor remainder of the logarithm and the remainder of \eqref{eq:local-sixth} divided by the Gaussian factor are $O(n^{-3}(1+L_n^{12}+e^{L_n^2/2}))=o(n^{-2})$, and so are their second differences. Using $\He_4''(z)=12(z^2-1)$, we obtain
\begin{equation}\label{eq:curvature-second}
 \Delta^2\log p_{n,k}
 =-\sigma_n^{-2}
  +\frac{\kappa_{4,n}}{2\sigma_n^6}(z_{n,k}^2-1)
  +o(n^{-2}).
\end{equation}
Since $\mu_n=n\lambda'(0)+\psi'(0)+O(\varepsilon_n)$ and $d_n=2\mu_n=n+O(1)$, we have $\lambda'(0)=1/2$ and $d_n=n+2\psi'(0)+O(\varepsilon_n)$. By \eqref{eq:symmetric-cumulants},
\begin{equation}\label{eq:ulc-cumulants}
 \sigma_n^2=vn+\psi''(0)+O(\varepsilon_n),\qquad
 \kappa_{4,n}=n\lambda^{(4)}(0)+O(1).
\end{equation}
Expanding the binomial ratio one order further, with $k-d_n/2=\sigma_n z_{n,k}$, gives
\begin{equation}\label{eq:binomial-first}
 \Delta^2\log\binom{d_n}k
 =-\frac4n+\frac{4(2\psi'(0)+1)-16vz_{n,k}^2}{n^2}
       +O\!\left(n^{-3}(1+L_n^4)+\varepsilon_n n^{-2}\right).
\end{equation}
Substituting \eqref{eq:ulc-cumulants} into \eqref{eq:curvature-second} and subtracting $r_*$ times \eqref{eq:binomial-first}, we find
\begin{equation}\label{eq:ulc-critical-curvature}
 \Delta^2\log\frac{p_{n,k}}{\binom{d_n}k^{r_*}}=\frac1{n^2}\biggl[\frac{\psi''(0)-v(2\psi'(0)+1)}{v^2}+4z_{n,k}^2+\frac{\lambda^{(4)}(0)(z_{n,k}^2-1)}{2v^3}\biggr]+o(n^{-2}),
\end{equation}
where the substitutions cost $O(n^{-3}(1+L_n^2)+\varepsilon_nn^{-2})=o(n^{-2})$. Under \eqref{eq:ulc-critical-condition}, the constant term of the bracket is negative and the coefficient of $z_{n,k}^2$ is nonpositive, so the bracket is bounded above by a negative constant on the whole range. This proves (2).
\end{proof}

\begin{theorem}\label{thm:heavy-ulc}
Fix $g\ge1$ and $0<c_0<\sqrt2$. The following statements hold uniformly for all $k$ such that $|k+j-d_n/2|\le c_0\,\sigma^H_{g,n}\sqrt{\log n}$ for $j=-1,0,1$.
\begin{enumerate}
\item For each fixed $0\le r<r_H$, the sequence $(b_{2k}(\Mbar_{g,n}))_{k=0}^{d_n}$ is strictly $r$-ultra-log-concave at $k$ for all sufficiently large $n$. For each fixed $r>r_H$, the strict reverse of \eqref{eq:ulc-log-criterion} holds at $k$ for all sufficiently large $n$.
\item At $r=r_H$, the sequence $(b_{2k}(\Mbar_{g,n}))_{k=0}^{d_n}$ is strictly $r_H$-ultra-log-concave at $k$ for all sufficiently large $n$.
\end{enumerate}
\end{theorem}

\begin{proof}
The hypotheses of Theorem~\ref{thm:symmetric-local} were verified in Section~\ref{subsec:betti-local}, with $v=\nu$ by \eqref{eq:lambda-psi-derivatives}. Moreover $d_n=n+O(1)$ and the mean is $d_n/2$. Lemma~\ref{lem:central-ulc}(1) with $r_*=1/(4\nu)=r_H$ proves (1). For (2), a direct computation from \eqref{eq:lambda-functions} and \eqref{eq:lambda-psi-derivatives} gives
\[
 \begin{aligned}
 \lambda_H^{(4)}(0)+8\nu^3&=-\frac{e(4e^2+9e-54)}{1080(e-2)^3}=-0.000141\ldots,\\
 \psi_H''(0)-\nu\bigl(2\psi_H'(0)+1\bigr)-\frac{\lambda_H^{(4)}(0)}{2\nu}&=-\frac{10(g-1)(7-2e)(3-e)+32e^2-219e+360}{120(e-2)(3-e)}.
 \end{aligned}
\]
Since $4e^2+9e-54=0.0207\ldots>0$, $7-2e>0$, $3-e>0$ and $32e^2-219e+360=1.14\ldots>0$, both expressions are negative for every $g\ge1$. Thus \eqref{eq:ulc-critical-condition} holds, and Lemma~\ref{lem:central-ulc}(2) proves (2).
\end{proof}

\begin{theorem}\label{thm:light-ulc}
Fix $g\ge1$ and $0<c_0<\sqrt2$. The following statements hold uniformly for all $k$ such that $|k+j-d_n/2|\le c_0\,\sigma^L_{g,n}\sqrt{\log n}$ for $j=-1,0,1$.
\begin{enumerate}
\item For each fixed $0\le r<3$, the sequence $(b_{2k}(\Mbar_{g,0|n}))_{k=0}^{d_n}$ is strictly $r$-ultra-log-concave at $k$ for all sufficiently large $n$. For each fixed $r>3$, the strict reverse of \eqref{eq:ulc-log-criterion} holds at $k$ for all sufficiently large $n$.
\item At $r=3$, the sequence $(b_{2k}(\Mbar_{g,0|n}))_{k=0}^{d_n}$ is strictly $3$-ultra-log-concave at $k$ for all sufficiently large $n$.
\end{enumerate}
\end{theorem}

\begin{proof}
The proof of (1) is the same as in the proof of Theorem~\ref{thm:heavy-ulc}, now with $v=1/12$ and $r_*=3$. Since $\lambda_L(s)=\log((e^s-1)/s)=\frac s2+\frac{s^2}{24}-\frac{s^4}{2880}+O(s^6)$, we have $\lambda_L^{(4)}(0)=-1/120$, and \eqref{eq:lambda-psi-derivatives} gives
\[
 \lambda_L^{(4)}(0)+8v^3=-\frac1{270},\qquad
 \psi_L''(0)-v\bigl(2\psi_L'(0)+1\bigr)-\frac{\lambda_L^{(4)}(0)}{2v}
 =-\frac{10g-9}{30}.
\]
Both are negative for $g\ge1$, and Lemma~\ref{lem:central-ulc} applies.
\end{proof}

\begin{example}[Losev--Manin spaces]\label{ex:losev-manin}
The critical case in Lemma~\ref{lem:central-ulc}(2) is not automatic. Let $Y_n$ have the normalized even Betti distribution of the Losev--Manin space $\Mbar_{0,2|n}$, of dimension $d_n=n-1$. Its even Betti polynomial is the Eulerian polynomial \cite[Theorem~2.3]{LosevManin}, and $\mu_n=(n-1)/2$. The generating function $\mathcal C$ of \eqref{eq:node-factors} has a simple pole at $w=\tau(t)$ with principal part $[t\tau(t)]^{-1}(1-w/\tau(t))^{-1}$, and its other poles tend to infinity as $t\to1$. Uniform coefficient extraction therefore gives \eqref{eq:symmetric-quasi} with an exponentially small error and
\[
 \lambda(s)=\log\frac{e^s-1}{s},\qquad \psi(s)=\lambda(s)-s.
\]
For $t=e^{i\theta}$ with $\delta\le|\theta|\le\pi$, the poles have modulus at least $|\theta|/(2|\sin(\theta/2)|)>1$, uniformly, so a Cauchy contour of radius greater than $1$ gives \eqref{eq:symmetric-off-zero}, the total Betti number being $n!$. Thus Theorem~\ref{thm:symmetric-local} applies, with
\[
 v=\frac1{12},\qquad \psi'(0)=-\frac12,\qquad
 \psi''(0)=\frac1{12},\qquad \lambda^{(4)}(0)=-\frac1{120},
\]
and Lemma~\ref{lem:central-ulc}(1) gives the threshold $r_*=3$. At $r=3$, however, the computation in the proof of Lemma~\ref{lem:central-ulc} gives, uniformly on every fixed bounded range of $z_{n,k}=(k-(n-1)/2)/\sigma_n$,
\[
 \Delta^2\log\frac{b_{2k}(\Mbar_{0,2|n})}{\binom{n-1}k^3}
 =\frac{16(6-z_{n,k}^2)}{5n^2}+O(n^{-3}),
\]
which is positive near the middle. Hence central $3$-ultra-log-concavity fails for $\Mbar_{0,2|n}$, although it holds for $\Mbar_{g,0|n}$ by Theorem~\ref{thm:light-ulc}. The difference comes from the constant terms $\psi'(0)$ and $\psi''(0)$, which are determined by the amplitude.
\end{example}

\section{Hodge types off the diagonal}\label{sec:hodge-refined}

Throughout this section, let $g\ge2$. The leading terms of the Hodge series \eqref{eq:hodge-pole} depend only on $uv$, so to leading order all cohomology lies on the Hodge diagonal. We now identify the next-order terms, which carry the off-diagonal classes. We first refine the generating series so that they record the Hodge difference (Section~\ref{subsec:refined-series}), then compute the contribution of an elliptic vertex (Section~\ref{subsec:elliptic-vertex}), and then sort the strata of $\Mbar_g$ by dimension (Section~\ref{subsec:small-strata}). The main result, Theorem~\ref{thm:hodge-asymmetry}, is proved in Section~\ref{subsec:hodge-asymmetry}. The results of Section~\ref{subsec:small-strata} are also used for the lower bound in Section~\ref{sec:tautological}.

\subsection{The Hodge-refined series}\label{subsec:refined-series}
Substitute $(u,v)=(-s\xi,-s\xi^{-1})$ in the full Hodge series:
\begin{equation}\label{eq:substituted-hodge-series}
 \widetilde{\mathscr L}_g(w;s,\xi)
 :=\mathscr L_g(w;-s\xi,-s\xi^{-1}),
 \qquad
 \widetilde{\mathscr H}_g(z;s,\xi)
 :=\mathscr H_g(z;-s\xi,-s\xi^{-1}).
\end{equation}
By purity,
\[
 \widetilde{\mathscr L}_g(w;s,\xi)
 =\sum_{n\ge0}\frac{w^n}{n!}
   \sum_{p,q}h^{p,q}(\Mbar_{g,0|n})\,s^{p+q}\xi^{p-q},
\]
and similarly for $\widetilde{\mathscr H}_g$, so $s$ records the degree and $\xi$ the Hodge difference. Since $uv=s^2$, Proposition~\ref{prop:substitution} gives
\[
 \widetilde{\mathscr H}_g(z;s,\xi)=\widetilde{\mathscr L}_g(T(z,s^2);s,\xi).
\]
For a function $F$ of $\xi$, let
\[
 F^{\mathrm{odd}}(\xi)=\frac{F(\xi)-F(-\xi)}2,
 \qquad
 F^{\mathrm{even}}(\xi)=\frac{F(\xi)+F(-\xi)}2.
\]
Since $p+q$ and $p-q$ have the same parity, these operations select the classes of odd and even degree, respectively, while the coefficient $[\xi^0]$ selects the Hodge diagonal.

For $Y=\Mbar_{g,n}$ or $Y=\Mbar_{g,0|n}$, choose a Hodge type $(p,q)$ with probability
\begin{equation}\label{eq:hodge-probability}
 \frac{h^{p,q}(Y)}{P_Y(1)}
 =\frac{h^{p,q}(Y)}{\sum_j b_j(Y)},
\end{equation}
and let $D^H_{g,n}$, respectively $D^L_{g,n}$, be the resulting Hodge difference $p-q$. Its joint generating function with the degree is
\[
 \sum_{p,q}\frac{h^{p,q}(Y)}{P_Y(1)}\,s^{p+q}\xi^{p-q}
 =\frac{E(Y;-s\xi,-s\xi^{-1})}{P_Y(1)}.
\]
By Hodge symmetry the distribution of the Hodge difference is symmetric about $0$. Every class of odd degree lies off the diagonal, and
\begin{equation}\label{eq:off-diagonal-odd-even}
 \sum_{p\ne q}h^{p,q}(Y)
 =\sum_k b_{2k+1}(Y)
  +\sum_{\substack{p\ne q\\p+q\ \mathrm{even}}}h^{p,q}(Y).
\end{equation}
Proposition~\ref{prop:hodge-pole} already shows that off-diagonal classes are rare. The leading term of $\widetilde{\mathscr L}_g(w;1,\xi)$ does not depend on $\xi$, so after subtracting the coefficient of $\xi^0$, what remains has denominator $(1-w)^{3g-4}$ and an entire numerator, uniformly on each annulus $r^{-1}\le|\xi|\le r$. The argument of Theorem~\ref{thm:coefficients}, followed by division by the Betti totals as in Corollary~\ref{cor:full-degree}, gives
\begin{equation}\label{eq:hodge-concentration-bound}
 \mathbb P(D^H_{g,n}\ne0)=O(n^{-1/2}),
 \qquad
 \mathbb P(D^L_{g,n}\ne0)=O(n^{-1}).
\end{equation}
To find the exact orders and the conditional distribution of the Hodge difference, we must compute the leading off-diagonal contribution.

\subsection{The elliptic vertex}\label{subsec:elliptic-vertex}
Off-diagonal classes come from vertices of positive genus. To determine their leading contribution, we first compute the vertex series in genus one. For $d=1,2$, consider smooth genus-one curves with $d$ ordered, pairwise distinct attaching points and $r$ labeled light markings that avoid the attaching points. Let $A_{1,d}(w;s,\xi)$ be the exponential generating series in $r$ of the Hodge--Deligne polynomials of these moduli spaces, evaluated at $(u,v)=(-s\xi,-s\xi^{-1})$. Equivalently, these spaces are the $r$-fold fiber powers over $\mathcal M_{1,d}$ of the universal elliptic curve with its $d$ attaching sections removed.

For the identity automorphism, the vertex series $\mathcal B_{\Gamma,1}$ in \eqref{eq:vertex-series} is the product of the individual vertex series. Thus $A_{1,d}$ is the factor associated with a genus-one vertex of valency $d$, after the substitution $(u,v)=(-s\xi,-s\xi^{-1})$. For example, if $\Gamma$ has one genus-one vertex of valency one and all its other vertices are rational and trivalent, then
\[
 \mathcal B_{\Gamma,1}(w;-s\xi,-s\xi^{-1})=A_{1,1}(w;s,\xi)e^{(2g-3)(s^2-2)w},
\]
since each of the $2g-3$ rational trivalent vertices contributes $e^{(s^2-2)w}$.

Recall from the proof of Lemma~\ref{lem:light-genus-one} that $I_n(u,v)$ is the Hodge--Deligne polynomial of the moduli space of smooth genus-one curves with $n$ labeled light markings. Taking the first marking as the origin identifies this space with the $(n-1)$-fold fiber power of the universal elliptic curve over $\mathcal M_{1,1}$, with no sections removed.

For $d=1$, we start with $I_{r+1}$ and exclude the origin from the $r$ light markings. Fixing $k$ specified light markings at the origin leaves the contribution $I_{r-k+1}$. Inclusion--exclusion gives
\[
 r!\,[w^r]A_{1,1}(w;s,\xi)
 =\sum_{k=0}^{r}(-1)^k\binom{r}{k}
 I_{r-k+1}(-s\xi,-s\xi^{-1}).
\]
Passing to exponential generating series yields the factor $\sum_{k\ge0}(-w)^k/k!=e^{-w}$. For $d=2$, first allow the second attaching point to move freely along with the $r$ light markings, giving $I_{r+2}$. The locus where the second attaching point equals the origin contributes $I_{r+1}$, so removing it gives $I_{r+2}-I_{r+1}$. Each light marking must then avoid the two distinct attaching points. Since each excluded marking has two possible locations, inclusion--exclusion gives the factor $e^{-2w}$. Thus
\begin{equation}\label{eq:punctured-elliptic-series}
 \begin{aligned}
 A_{1,1}(w;s,\xi)
 &=e^{-w}\sum_{r\ge0}I_{r+1}(-s\xi,-s\xi^{-1})\frac{w^r}{r!},\\
 A_{1,2}(w;s,\xi)
 &=e^{-2w}\sum_{r\ge0}
   \bigl(I_{r+2}(-s\xi,-s\xi^{-1})-I_{r+1}(-s\xi,-s\xi^{-1})\bigr)
   \frac{w^r}{r!}.
 \end{aligned}
\end{equation}
Since $I_n^{\mathrm{even}}$ depends only on $uv=s^2$, the even parts of both series are independent of $\xi$. These series are entire in $w$, locally uniformly for $(s,\xi)\in\mathbb C\times\mathbb C^*$, by the bounds in the proof of Lemma~\ref{lem:light-genus-one}.

The odd part of $A_{1,1}$ has a closed form in terms of the modified Bessel functions
\[
 \mathbb{I}_m(2w)=\sum_{j\ge0}\frac{w^{m+2j}}{j!\,(m+j)!},\qquad m\ge0,
\]
and the cusp-form dimensions $\beta_k=\dim S_k(\mathrm{SL}_2(\mathbb Z))$ for even $k\ge4$. Define
\begin{equation}\label{eq:cusp-bessel-series}
 \mathcal E(w;s,\xi)=\sum_{\substack{m\ge3\\m\ \mathrm{odd}}}\beta_{m+1}\,\frac{m}{w}\,\mathbb{I}_m(2sw)\,(\xi^m+\xi^{-m}),
\end{equation}
with the value at $w=0$ defined by continuity.

\begin{lemma}\label{lem:elliptic-odd}
The series $\mathcal E$ is entire in $(w,s)$ and holomorphic in $\xi\ne0$. With $t=s^2$ and $\alpha_{ij}=\binom{2j}{i}-\binom{2j}{i-1}$,
\[
 A_{1,1}^{\mathrm{odd}}(w;s,\xi)=e^{tw}\sum_{j\ge1}\frac{s^{2j+1}w^{2j}}{(2j)!}\sum_{i=0}^{j-1}\alpha_{ij}\beta_{2j-2i+2}\bigl(\xi^{2j-2i+1}+\xi^{-(2j-2i+1)}\bigr)=e^{tw}\,\mathcal E(w;s,\xi).
\]
In particular, $A_{1,1}^{\mathrm{odd}}$ contains only odd powers of $\xi$. At $s=1$ and real $w>0$ its Laurent coefficients are nonnegative, and
\begin{equation}\label{eq:double-elliptic-positivity}
 [\xi^{11}]A_{1,1}^{\mathrm{odd}}(w;1,\xi)\ge e^w\frac{w^{10}}{10!},\qquad [\xi^{22}]\bigl(A_{1,1}^{\mathrm{odd}}(1;1,\xi)\bigr)^2\ge\frac{e^2}{(10!)^2}>0.
\end{equation}
\end{lemma}

\begin{proof}
At $(u,v)=(-s\xi,-s\xi^{-1})$ we have $uv=t$ and $u^m+v^m=-s^m(\xi^m+\xi^{-m})$ for odd $m$. Inserting this into the formula for $I_n^{\mathrm{odd}}$ and summing
\[
 \sum_{r\ge2j}\binom{r}{2j}(1+t)^{r-2j}\frac{w^r}{r!}=e^{(1+t)w}\frac{w^{2j}}{(2j)!}
\]
gives the first expression, the factor $e^{-w}$ in \eqref{eq:punctured-elliptic-series} turning $e^{(1+t)w}$ into $e^{tw}$.

For the second, fix an odd $m\ge3$ and collect the terms with $2j-2i+1=m$, that is, $i=j-(m-1)/2$. Writing $j=(m-1)/2+l$, the definition of $\alpha_{ij}$ gives
\[
 \alpha_{l,(m-1)/2+l}=\binom{m-1+2l}{l}-\binom{m-1+2l}{l-1}=\frac{m\,(m-1+2l)!}{l!\,(m+l)!},
\]
so that
\[
 \sum_{j\ge(m-1)/2}\alpha_{j-(m-1)/2,j}\frac{s^{2j+1}w^{2j}}{(2j)!}=\sum_{l\ge0}\frac{m\,s^{m+2l}w^{m-1+2l}}{l!\,(m+l)!}=\frac{m}{w}\,\mathbb{I}_m(2sw).
\]

For $|w|\le R$ and $|s|\le S$ with $R,S\ge1$, the Bessel series gives
\[
 \left|\frac{m}{w}\mathbb{I}_m(2sw)\right|\le\frac{mS^mR^{m-1}}{m!}e^{S^2R^2}.
\]
Since $\beta_{m+1}=O(m)$, the series \eqref{eq:cusp-bessel-series} converges normally on compact subsets of $\mathbb C^2\times\mathbb C^*$. Hence $\mathcal E$ is entire in $(w,s)$ and holomorphic in $\xi\ne0$. These bounds also justify the rearrangements above.

Since $\alpha_{ij}\ge0$ for $i\le j$ and $\beta_k\ge0$, all Laurent coefficients are nonnegative at $s=1$ and $w>0$. The term $i=0$, $j=5$ has $\alpha_{0,5}=\beta_{12}=1$ and gives \eqref{eq:double-elliptic-positivity}.
\end{proof}

\subsection{Strata of small dimension}\label{subsec:small-strata}
Recall from Proposition~\ref{prop:hodge-pole} that the full light Hodge series has a pole of order $3g-3$. We next identify the stable cores that can contribute off-diagonal terms at pole orders $3g-4$ and $3g-5$. We also track the even off-diagonal part and, in particular, the coefficient of $\xi^{22}$, whose positive leading term will give a lower bound for the even off-diagonal contribution. The following lemma classifies the relevant strata of dimensions one and two and computes their contributions.

\begin{lemma}\label{lem:small-strata}
Let $\Gamma$ be a connected stable graph of genus $g\ge2$ without legs, and let $\sigma\in\Aut\Gamma$. In (3) and (4), we evaluate at $(u,v)=(-s\xi,-s\xi^{-1})$, put $t=s^2$, and suppress the arguments of $A_{1,d}(w;s,\xi)$, $\mathcal E(w;s,\xi)$, and $\mathcal C(w,t)$. The contributions in the tables omit the graph weight $1/|\Aut\Gamma|$.
\begin{enumerate}
\item The term of $(\Gamma,\sigma)$ in \eqref{eq:graph-sum} is an entire function divided by $(1-W)^{p(\sigma)}$, where
\[
 p(\sigma)\le |E(\Gamma)|=3g-3-\dim\mathcal M_\Gamma,
\]
with equality only if $\sigma=1$.

\item If every vertex of $\Gamma$ has genus zero, the term depends on $(u,v)$ only through $uv$. After the substitution \eqref{eq:substituted-hodge-series}, it is independent of $\xi$.

\item If $\Gamma$ has a vertex of positive genus and $\dim\mathcal M_\Gamma=1$, then it has exactly one genus-one vertex of valency one, and all other vertices are rational and trivalent. We call such a $\Gamma$ an \emph{elliptic-leaf graph}. It has $3g-4$ edges and $2g-3$ rational vertices. Its identity contribution is as follows.
\begin{center}
\renewcommand{\arraystretch}{1.2}
\begin{tabular}{ll}
\hline
Part & Identity contribution \\
\hline
Full & $e^{(2g-3)(t-2)w}A_{1,1}\mathcal C^{3g-4}$ \\
Odd & $e^{2((g-1)(t-2)+1)w}\mathcal E\,\mathcal C^{3g-4}$ \\
Even off-diagonal & $0$ \\
\hline
\end{tabular}
\end{center}
Moreover, for every $\sigma\in\Aut\Gamma$, the corresponding term after the substitution \eqref{eq:substituted-hodge-series} has the form $A_{1,1}(w;s,\xi)F_\sigma(w;s)$, with $F_\sigma$ independent of $\xi$. Thus its even part is independent of $\xi$.

\item If $\Gamma$ has a vertex of positive genus and $\dim\mathcal M_\Gamma=2$, then the possibilities are listed below. Every vertex not listed in the first column is rational and trivalent. In all three cases, $\Gamma$ has $3g-5$ edges and $2g-4$ rational vertices. Its identity contribution is
\[
 e^{(2g-4)(t-2)w}\mathcal C^{3g-5}
\]
times the factor in the second column.
\begin{center}
\renewcommand{\arraystretch}{1.2}
\begin{tabular}{lll}
\hline
Vertex types other than $(0,3)$ & Factor & $\xi^{22}$-coefficient of the factor \\
\hline
$(1,1),(1,1)$ & $A_{1,1}^{\,2}$ & $[\xi^{22}](A_{1,1}^{\mathrm{odd}})^2$ \\
$(1,2)$ & $A_{1,2}$ & $0$ \\
$(1,1),(0,4)$ ($g\ge3$) & $(t-2)e^{-w}A_{1,1}$ & $0$ \\
\hline
\end{tabular}
\end{center}
Here the vertex types are the pairs $(g_v,d_v)$. We call the graphs in the first row \emph{double-leaf graph}. Only these graphs can contribute even powers $\xi^{2k}$ with $k\ne0$ to the identity term. In the third case, the $2g-4$ rational vertices consist of one vertex of valency four and $2g-5$ trivalent vertices.
\end{enumerate}
\end{lemma}

\begin{proof}
Assertion (1) follows from Lemmas~\ref{lem:node-factors} and~\ref{lem:entire-vertices} and \eqref{eq:stratum-dimension}. If $p(\sigma)=|E(\Gamma)|$, then $\sigma$ fixes every half-edge, hence every vertex, and $\sigma=1$.

For (2), stratifying the configurations of markings on a rational vertex by their coincidence partitions expresses the vertex spaces in terms of products of moduli spaces $\mathcal M_{0,k}$. Their compactly supported cohomology is Hodge--Tate, so the traces $E_\sigma$ depend only on $uv$. The chain factors also depend only on $uv$, by Lemma~\ref{lem:node-factors}.

For the classifications in (3) and (4), use
\[
 \dim\mathcal M_\Gamma=\sum_v(3g_v-3+d_v).
\]
A rational vertex contributes $d_v-3\ge0$, a genus-one vertex contributes $d_v\ge1$, and a vertex of genus at least two contributes at least $3$. This gives the listed possibilities. The edge counts follow from \eqref{eq:stratum-dimension}, and the rational vertex counts follow by summing the valencies.

For the identity contributions, the vertex series multiply and each edge contributes $\mathcal C$. A rational trivalent vertex contributes $e^{(t-2)w}$. For a rational vertex of valency four, taking the first three attaching points to be $0,1,\infty$ identifies the unpunctured $r$-fold family with $\mathcal M_{0,4}\times(\mathbb P^1)^r$. Removing the four attaching sections from each coordinate by inclusion--exclusion gives the Hodge--Deligne polynomial $(t-2)(t-3)^r$. Its exponential generating series is therefore
\[
 (t-2)e^{(t-3)w}.
\]
Multiplying these rational factors by the elliptic factors $A_{1,1}$ or $A_{1,2}$ gives the full contributions in both tables. In the third row of (4), there are $2g-5$ rational trivalent vertices, and
\[
 e^{(2g-5)(t-2)w}(t-2)e^{(t-3)w}=e^{(2g-4)(t-2)w}(t-2)e^{-w}.
\]
The tables omit only the graph weight from \eqref{eq:graph-sum}. The cohomology defining $A_{1,d}$ already accounts for the stabilizers of the pointed elliptic curves.

By Lemma~\ref{lem:elliptic-odd}, $A_{1,1}^{\mathrm{odd}}=e^{tw}\mathcal E$. Together with
\[
 (2g-3)(t-2)+t=2((g-1)(t-2)+1),
\]
this gives the odd contribution in (3). By \eqref{eq:punctured-elliptic-series}, the even parts of $A_{1,1}$ and $A_{1,2}$ are independent of $\xi$, which proves the vanishing assertions in both tables. In the first row of (4), $(A_{1,1}^{\mathrm{even}})^2$ is independent of $\xi$, while $A_{1,1}^{\mathrm{even}}A_{1,1}^{\mathrm{odd}}$ contains only odd powers of $\xi$. Thus every nonconstant even Laurent coefficient comes from $(A_{1,1}^{\mathrm{odd}})^2$.

Finally, every automorphism fixes the unique elliptic leaf and its attaching flag, so $A_{1,1}$ factors out of every term of an elliptic-leaf graph. The remaining rational vertex and chain factors are independent of $\xi$, by (2) and Lemma~\ref{lem:node-factors}. This proves the assertion in (3) for arbitrary $\sigma$, including the independence of its even part from $\xi$.
\end{proof}

\begin{lemma}\label{lem:leaf-weight}
The sum of $1/|\Aut\Gamma|$ over the elliptic-leaf graphs of genus $g$ is
\begin{equation}\label{eq:elliptic-leaf-weight}
 \ell_2=\frac12,\qquad \ell_g=(3g-6)\,c_{g-1}\quad(g\ge3),
\end{equation}
where $c_g$ is defined in \eqref{eq:core-constant}.
\end{lemma}

\begin{proof}
For $g=2$, there is one elliptic-leaf graph. The leaf is attached to a rational vertex carrying a loop, and reversing the loop gives $|\Aut\Gamma|=2$. For $g\ge3$, remove the leaf and its edge and suppress the resulting bivalent rational vertex. That vertex cannot carry a loop, since otherwise the graph would be the genus-two graph above, so the result is a trivalent graph of genus $g-1$ with a distinguished unoriented edge. Conversely, subdividing a distinguished edge of a trivalent graph of genus $g-1$ and attaching an elliptic leaf at the new vertex recovers the original graph, with the same automorphism group. Summing over the orbits of edges gives $\ell_g=\sum_{\Gamma\in\mathscr G_{g-1}}|E(\Gamma)|/|\Aut\Gamma|=(3g-6)c_{g-1}$.
\end{proof}

\begin{proposition}\label{prop:next-poles}
With $t=s^2$, there are functions $\mathscr K_g^{\mathrm{odd}}$ and $\mathscr K_g^{\mathrm{even}}$, entire in $w$ and holomorphic in $(s,\xi)\in\mathbb C\times\mathbb C^*$, locally uniformly, such that
\begin{equation}\label{eq:odd-hodge-pole}
 \widetilde{\mathscr L}_g^{\mathrm{odd}}(w;s,\xi)
 =
 \frac{\ell_g\, e^{2((g-1)(t-2)+1)w}\,\mathcal E(w;s,\xi)}
      {(1-W(w,t))^{3g-4}}
 +\frac{\mathscr K_g^{\mathrm{odd}}(w;s,\xi)}
       {(1-W(w,t))^{3g-5}}
\end{equation}
and
\begin{equation}\label{eq:even-off-diagonal-pole}
 \widetilde{\mathscr L}_g^{\mathrm{even}}(w;s,\xi)
 -[\xi^0]\widetilde{\mathscr L}_g^{\mathrm{even}}(w;s,\xi)
 =\frac{\mathscr K_g^{\mathrm{even}}(w;s,\xi)}
       {(1-W(w,t))^{3g-5}}.
\end{equation}
Moreover, there are $a_g>0$ and an entire function $\mathcal N_g$ such that
\begin{equation}\label{eq:fixed-hodge-difference-pole}
 [\xi^{22}]\widetilde{\mathscr L}_g(w;1,\xi)
 =\frac{a_g}{(1-w)^{3g-5}}
  +\frac{\mathcal N_g(w)}{(1-w)^{3g-6}}.
\end{equation}
\end{proposition}

\begin{proof}
By Lemma~\ref{lem:small-strata} (2), rational cores contribute neither to the odd part nor to the off-diagonal part. Among the cores with a vertex of positive genus, Lemma~\ref{lem:small-strata} (1) and (3) shows that only the identity terms of elliptic-leaf graphs can contribute at pole order $3g-4$, while all other terms have pole order at most $3g-5$. Summing the odd contributions in the table of Lemma~\ref{lem:small-strata} (3), with the graph weights computed in Lemma~\ref{lem:leaf-weight}, gives the leading term of \eqref{eq:odd-hodge-pole}. Putting the remaining odd terms over the common denominator $(1-W)^{3g-5}$ gives the stated entire remainder.

By Lemma~\ref{lem:small-strata} (3), every term of an elliptic-leaf graph, including the nonidentity terms, has an even part independent of $\xi$. Thus the even off-diagonal contribution comes only from strata of dimension at least two, whose terms have pole order at most $3g-5$ by Lemma~\ref{lem:small-strata} (1). Putting these terms over the common denominator $(1-W)^{3g-5}$ proves \eqref{eq:even-off-diagonal-pole}.

For \eqref{eq:fixed-hodge-difference-pole}, take $s=1$, so that $t=1$ and $\mathcal C=(1-w)^{-1}$, and extract the coefficient of $\xi^{22}$. Rational cores and elliptic-leaf graphs contribute nothing, by Lemma~\ref{lem:small-strata} (2) and (3). Terms from strata of dimension at least three, and nonidentity terms from strata of dimension two, have pole order at most $3g-6$, by Lemma~\ref{lem:small-strata} (1). The table in Lemma~\ref{lem:small-strata} (4) therefore shows that only double-leaf identity terms remain at order $3g-5$. Their leading coefficient at $w=1$ is
\[
 a_g=e^{-(2g-4)}[\xi^{22}]\bigl(A_{1,1}^{\mathrm{odd}}(1;1,\xi)\bigr)^2\sum_{\Gamma\text{ double-leaf}}\frac1{|\Aut\Gamma|},
\]
where the sum runs over the isomorphism classes of double-leaf graphs of genus $g$.

Double-leaf graphs exist for every $g\ge2$. For $g=2$, join two elliptic vertices by an edge. Subdividing an edge twice and joining the two new rational vertices by an additional edge raises the genus by one while preserving the required valencies. By \eqref{eq:double-elliptic-positivity}, the Laurent coefficient in the formula for $a_g$ is strictly positive, so $a_g>0$.

Expanding the double-leaf numerators at $w=1$ and putting the remaining terms over the denominator $(1-w)^{3g-6}$ gives \eqref{eq:fixed-hodge-difference-pole}. Laurent coefficient extraction is performed by integration over $|\xi|=1$, which preserves holomorphy. The differences between the double-leaf numerators and their values at $w=1$ are divisible by $1-w$ with entire quotients, so the resulting remainder numerator $\mathcal N_g$ is entire.
\end{proof}

\subsection{The Hodge difference among off-diagonal classes}\label{subsec:hodge-asymmetry}
Let
\begin{equation}\label{eq:hodge-conditional-law}
 Z_Q=\sum_{m\ge3,\ m\ \mathrm{odd}}\beta_{m+1}\,m\,\mathbb I_m(2),\qquad q_m=q_{-m}=\frac{\beta_{m+1}\,m\,\mathbb I_m(2)}{2Z_Q}\quad(m\ge3\text{ odd}),
\end{equation}
and $q_d=0$ for all other integers $d$. The series converges, and $Z_Q>0$ since $\beta_{12}=1$. Thus $(q_d)_{d\in\mathbb Z}$ defines a probability distribution $Q$ on $\mathbb Z$, whose Laurent generating function is
\[
 \widehat Q(\xi)=\sum_{d\in\mathbb Z}q_d\xi^d=\frac{\mathcal E(1;1,\xi)}{2Z_Q}.
\]
Finally, let
\begin{equation}\label{eq:hodge-rarity-constants}
 \gamma_g^H=\frac{2e^2\ell_gZ_Q}{c_g}\sqrt{\frac{2(e-2)}e}\,\frac{\Gamma((3g-3)/2)}{\Gamma((3g-4)/2)},\qquad \gamma_g^L=\frac{2e^2\ell_gZ_Q}{c_g}\,(3g-4).
\end{equation}

\begin{theorem}\label{thm:hodge-asymmetry}
Fix $g\ge2$.
\begin{enumerate}
\item The off-diagonal masses satisfy
\begin{equation}\label{eq:hodge-rarity}
 \mathbb P(D^H_{g,n}\ne0)=\gamma_g^H n^{-1/2}+O(n^{-1}),\qquad \mathbb P(D^L_{g,n}\ne0)=\gamma_g^L n^{-1}+O(n^{-2}).
\end{equation}
The same formulas hold for the probability that $p+q$ is odd, and
\begin{equation}\label{eq:even-off-diagonal-mass}
 \mathbb P(D^H_{g,n}\ne0,\ p+q\ \mathrm{even})=O(n^{-1}),\qquad \mathbb P(D^L_{g,n}\ne0,\ p+q\ \mathrm{even})=O(n^{-2}).
\end{equation}
\item For $J\in\{H,L\}$, the conditional distribution of $D^J_{g,n}$ given $D^J_{g,n}\ne0$ converges to $Q$. More precisely, for every fixed $r>1$,
\begin{equation}\label{eq:hodge-conditional-rate}
 \sum_{d\ne0}r^{|d|}\left|\mathbb P(D^J_{g,n}=d\mid D^J_{g,n}\ne0)-q_d\right|=\begin{cases}O_r(n^{-1/2}),&J=H,\\O_r(n^{-1}),&J=L.\end{cases}
\end{equation}
Moreover, for all sufficiently large $n$ and all integers $k\ge1$,
\begin{equation}\label{eq:hodge-exponential-tails}
 \mathbb P(|D^H_{g,n}|\ge k)\le C_rn^{-1/2}r^{-k},\qquad \mathbb P(|D^L_{g,n}|\ge k)\le C_rn^{-1}r^{-k}.
\end{equation}
\item For every $x\in\mathbb R$ and $d\in\mathbb Z$,
\begin{equation}\label{eq:hodge-conditional-joint}
 \mathbb P\!\left(\frac{p+q-d_n}{\sqrt{4\nu n}}\le x,\ D^H_{g,n}=d\ \middle|\ D^H_{g,n}\ne0\right)\longrightarrow\Phi(x)q_d,
\end{equation}
where $\Phi$ is the standard Gaussian distribution function. The same limit holds in the light case with $\sqrt{4\nu n}$ replaced by $\sqrt{n/3}$ and $D^H_{g,n}$ by $D^L_{g,n}$. Thus, conditionally on being off the Hodge diagonal, the normalized degree and the Hodge difference are asymptotically independent.
\end{enumerate}
\end{theorem}

\begin{proof}
\emph{Coefficient asymptotics at $s=1$.} At $s=1$, replacing the numerator $\ell_g e^{-(2g-4)w}\mathcal E(w;1,\xi)$ in \eqref{eq:odd-hodge-pole} by its value $\ell_g e^{-(2g-4)}\mathcal E(1;1,\xi)$ at $w=1$ changes only the remainder with denominator $(1-w)^{3g-5}$. The transfer theorem \cite[Theorem~VI.4]{FS}, together with the uniformity lemma \cite[Lemma~IX.2]{FS}, gives, uniformly on each annulus $r^{-1}\le|\xi|\le r$,
\begin{equation}\label{eq:odd-light-coefficients}
 [w^n]\widetilde{\mathscr L}_g^{\mathrm{odd}}(w;1,\xi)=\frac{\ell_g e^{-(2g-4)}}{\Gamma(3g-4)}\,\mathcal E(1;1,\xi)\,n^{3g-5}+O_r(n^{3g-6}).
\end{equation}
For the heavy series, Lemma~\ref{lem:inverse} gives
\[
 1-T(z,1)=\sqrt{\frac{2(e-2)}e}\left(1-\frac z{e-2}\right)^{1/2}+O\!\left(1-\frac z{e-2}\right),
\]
and by \eqref{eq:matching} the function $1-T(z,1)$ has no zeros in the $\Delta$-domain of Lemma~\ref{lem:inverse}. Composing \eqref{eq:odd-hodge-pole} with $T(z,1)$ and applying the transfer theorem \cite[Theorem~VI.4]{FS}, with uniformity as in \cite[Lemma~IX.2]{FS}, gives
\begin{equation}\label{eq:odd-heavy-coefficients}
 [z^n]\widetilde{\mathscr H}_g^{\mathrm{odd}}(z;1,\xi)=(e-2)^{-n}\biggl(\frac{\ell_g e^{-(2g-4)}}{\Gamma((3g-4)/2)}\left(\frac{2(e-2)}e\right)^{-\frac{3g-4}{2}}\mathcal E(1;1,\xi)\,n^{\frac{3g-6}{2}}+O_r\bigl(n^{\frac{3g-7}{2}}\bigr)\biggr).
\end{equation}
The numerators are bounded on compact sets, and the continuation domain of $T(z,1)$ does not depend on $\xi$, so these estimates are uniform on the annulus. They are additive, so they remain valid at zeros of $\mathcal E(1;1,\xi)$. For $g=2$, the light remainder is a simple pole and the heavy remainder has exponent $-1/2$, and the same bounds hold.

\medskip

\emph{Proof of (1).} By Theorem~\ref{thm:coefficients} and \eqref{eq:full-even-ratios}, the total Betti numbers are $A_g(1)\,n!\,n^{3g-4}(1+O(n^{-1}))$ and $C_g(1)\,n!\,n^{(3g-5)/2}(e-2)^{-n}(1+O(n^{-1/2}))$ in the light and heavy cases, respectively. Multiplying \eqref{eq:odd-light-coefficients} and \eqref{eq:odd-heavy-coefficients} by $n!$ and dividing by the corresponding totals and substituting the amplitudes in \eqref{eq:amplitudes}, we obtain
\begin{equation}\label{eq:odd-hodge-probability-series}
\begin{aligned}
 \mathbb E\bigl[\xi^{D^H_{g,n}}\mathbf 1_{\{p+q\ \mathrm{odd}\}}\bigr]&=\gamma_g^H n^{-1/2}\widehat Q(\xi)+O_r(n^{-1}),\\
 \mathbb E\bigl[\xi^{D^L_{g,n}}\mathbf 1_{\{p+q\ \mathrm{odd}\}}\bigr]&=\gamma_g^L n^{-1}\widehat Q(\xi)+O_r(n^{-2}).
\end{aligned}
\end{equation}
Applying the same transfer to \eqref{eq:even-off-diagonal-pole} and dividing by the corresponding total Betti numbers shows that, for $J\in\{H,L\}$, the expectation $\mathbb E\bigl[\xi^{D^J_{g,n}}\mathbf 1_{\{p+q\ \mathrm{even},\,p\ne q\}}\bigr]$ is $O_r(n^{-1})$ in the heavy case and $O_r(n^{-2})$ in the light case. Since $p+q$ odd implies $p\ne q$, adding this even off-diagonal contribution shows that the corresponding identity in \eqref{eq:odd-hodge-probability-series} also holds with $\mathbf 1_{\{p+q\ \mathrm{odd}\}}$ replaced by $\mathbf 1_{\{D^J_{g,n}\ne0\}}$, with the same error order. Using $\widehat Q(1)=1$, evaluation at $\xi=1$ proves \eqref{eq:hodge-rarity} and \eqref{eq:even-off-diagonal-mass}. In particular, the conditioning events have positive probability for all sufficiently large $n$.

\emph{Proof of (2).} Dividing the off-diagonal Laurent generating functions obtained in (1) by their values at $\xi=1$ gives the conditional Laurent generating functions. These converge to $\widehat Q$, uniformly on each annulus $r^{-1}\le|\xi|\le r$, with errors $O_r(n^{-1/2})$ and $O_r(n^{-1})$ in the heavy and light cases, respectively. For each $d\ne0$, the coefficient of $\xi^d$ in the corresponding difference is $\mathbb P(D^J_{g,n}=d\mid D^J_{g,n}\ne0)-q_d$. Given $r>1$, choose $R>r$. Cauchy's estimates on $|\xi|=R$ and $|\xi|=R^{-1}$ bound these coefficient errors by $C_Rn^{-1/2}R^{-|d|}$ and $C_Rn^{-1}R^{-|d|}$, respectively. Multiplying by $r^{|d|}$ and summing the convergent geometric series $\sum_{d\ne0}(r/R)^{|d|}$ gives \eqref{eq:hodge-conditional-rate}. Finally, the off-diagonal Laurent generating functions have nonnegative coefficients, so evaluation at $\xi=r$ and $\xi=r^{-1}$ gives \eqref{eq:hodge-exponential-tails}.

\emph{Proof of (3).} Let $s$ vary in a complex neighborhood of $1$, and put $t=s^2$. The leading numerator of \eqref{eq:odd-hodge-pole} at $w=\tau(t)$ is $\ell_ge^{2((g-1)(t-2)+1)\tau(t)}\mathcal E(\tau(t);s,\xi)$. The leading local terms of the denominators are $t\tau(t)(1-w/\tau(t))$ in the light case and $(2\rho(t))^{1/2}e^{(t-2)\tau(t)/2}(1-z/\rho(t))^{1/2}$ in the heavy case, by \eqref{eq:composed-denominator}. The convergence of \eqref{eq:cusp-bessel-series}, the entire remainders of Proposition~\ref{prop:next-poles}, and the uniform domains and bounds from the proof of Theorem~\ref{thm:coefficients} give coefficient estimates uniformly in $s$ and on compact $\xi$-annuli. For each $n$, dividing the off-diagonal coefficient by their values at $(s,\xi)=(1,1)$ gives
\begin{equation}\label{eq:hodge-conditional-quasi-powers}
\begin{aligned}
 \mathbb E\bigl[s^{p+q}\xi^{D^H_{g,n}}\bigm|D^H_{g,n}\ne0\bigr]&=\left(\frac{e-2}{\rho(s^2)}\right)^n\bigl(\mathfrak a_g^H(s,\xi)+O_r(n^{-1/2})\bigr),\\
 \mathbb E\bigl[s^{p+q}\xi^{D^L_{g,n}}\bigm|D^L_{g,n}\ne0\bigr]&=\tau(s^2)^{-n}\bigl(\mathfrak a_g^L(s,\xi)+O_r(n^{-1})\bigr),
\end{aligned}
\end{equation}
where $\mathfrak a_g^H$ and $\mathfrak a_g^L$ are holomorphic for $s$ near $1$ and $\xi\ne0$, with $\mathfrak a_g^J(1,\xi)=\widehat Q(\xi)$ for $J\in\{H,L\}$. By the proof of Theorem~\ref{thm:clt},
\[
 \log\frac{e-2}{\rho(e^{2h})}=h+2\nu h^2+O(h^3),\qquad -\log\tau(e^{2h})=h+\frac{h^2}{6}+O(h^3).
\]
For fixed real $a,b$, let $\xi=e^{ib}$, and take $s=\exp(ia/\sqrt{4\nu n})$ in the heavy case and $s=\exp(ia/\sqrt{n/3})$ in the light case. Multiplying \eqref{eq:hodge-conditional-quasi-powers} by $s^{-d_n}$ and using $d_n=n+O(1)$ shows that the conditional joint characteristic functions converge to $e^{-a^2/2}\widehat Q(e^{ib})$. By L\'evy's continuity theorem, the conditional joint distributions converge to the product of the standard Gaussian distribution and $Q$. Since the Hodge difference is integer-valued and the Gaussian distribution function is continuous, applying this convergence to $(-\infty,x]\times(d-\tfrac12,d+\tfrac12)$ gives \eqref{eq:hodge-conditional-joint} for every $x\in\mathbb R$ and $d\in\mathbb Z$, and likewise in the light case.
\end{proof}

Theorem~\ref{thm:hodge-asymmetry} describes both the rarity and the Hodge distribution of off-diagonal classes. Their total mass tends to zero at rate $n^{-1/2}$ for $\Mbar_{g,n}$ and $n^{-1}$ for $\Mbar_{g,0|n}$, while their conditional Hodge difference has a universal limiting distribution $Q$, independent of the genus. On each line $p-q=d$ with $q_d>0$, the degree converges to the standard Gaussian after the normalization in (3), expressing the asymptotic independence of the normalized degree and the Hodge difference. The universality of $Q$ reflects the leading contribution of a single elliptic leaf, with the remaining Tate factors shifting $p$ and $q$ equally. The even-degree off-diagonal classes are of lower order, so the Hodge difference among odd-degree classes also converges in distribution to $Q$.

\begin{remark}[Interpretation of $Q$]\label{rem:skellam-interpretation}
Let $\operatorname{Sk}(j)=e^{-2}\mathbb I_{|j|}(2)$ denote the probability mass function of the difference of two independent Poisson variables of mean $1$. The Bessel recurrence $\mathbb I_{m-1}(2)-\mathbb I_{m+1}(2)=m\mathbb I_m(2)$ rewrites the definition of $q_m$ in \eqref{eq:hodge-conditional-law} as
\[
 q_m=q_{-m}=\frac{e^2\beta_{m+1}}{2Z_Q}\bigl(\operatorname{Sk}(m-1)-\operatorname{Sk}(m+1)\bigr),\qquad m\ge3\text{ odd}.
\]

The difference of Skellam probabilities in this formula has a representation-theoretic interpretation. Let $V$ be the standard representation of $\mathrm{SL}_2$, with character $\operatorname{ch}(V)(\xi)=\xi+\xi^{-1}$. Then
\[
 e^{-2}\sum_{k\ge0}\frac{\operatorname{ch}(V^{\otimes k})(\xi)}{k!}=\exp(\xi+\xi^{-1}-2)=\sum_{j\in\mathbb Z}\operatorname{Sk}(j)\xi^j.
\]
In an $\mathrm{SL}_2$-representation, subtracting the multiplicity of weight $m+1$ from that of weight $m-1$ gives the multiplicity of $\operatorname{Sym}^{m-1}V$, since the contributions from higher highest weights cancel. Thus $\operatorname{Sk}(m-1)-\operatorname{Sk}(m+1)$ is the weighted multiplicity of $\operatorname{Sym}^{m-1}V$ in the exponential sum above.

For odd $m\ge3$, Eichler--Shimura theory gives a cusp contribution of dimension $\beta_{m+1}$ for each Hodge difference $\pm m$ in the cohomology of the associated local system $\operatorname{Sym}^{m-1}V$ on $\mathcal M_{1,1}$. Multiplying the weighted highest-weight multiplicities by these cusp-form dimensions and normalizing therefore gives the probabilities $q_m=q_{-m}$.

Since the first nonzero cusp-form dimension is $\beta_{12}=1$, we have $q_d=0$ for $|d|<11$, and numerically
\[
q_{11}+q_{-11}=\frac{11\mathbb I_{11}(2)}{Z_Q}\approx0.99996.
\]
\end{remark}

\section{Tautological classes}\label{sec:tautological}

Let $RH^{2k}(\Mbar_{g,n})\subset H^{2k}(\Mbar_{g,n})$ denote the image of the tautological group $R^k(\Mbar_{g,n})$ under the cycle class map, and let
\[
 H^{\taut}_{g,n}(t)=\sum_k\dim_{\mathbb Q}RH^{2k}(\Mbar_{g,n})\,t^k.
\]
Let $X^{\taut}_{g,n}$ be the random variable with PGF $H^{\taut}_{g,n}(t)/H^{\taut}_{g,n}(1)$. In genus one, Petersen's theorem gives $H^{\taut}_{1,n}=H_{1,n}$ \cite{Petersen}. In this section we show that, for every $g\ge2$, the tautological classes still form a proportion $1-\Theta(n^{-1})$ of the even cohomology.

\begin{theorem}\label{thm:tautological-dominance}
Fix $g\ge2$. There are constants $c'_g,C'_g>0$ such that, for all sufficiently large $n$,
\begin{equation}\label{eq:tautological-fraction}
 \frac{c'_g}{n}
 \le 1-\frac{H^{\taut}_{g,n}(1)}{H_{g,n}(1)}
 \le\frac{C'_g}{n}.
\end{equation}
Consequently,
\begin{equation}\label{eq:tautological-tv}
 \frac12\sum_k
 \left|\mathbb P(X^{\taut}_{g,n}=k)
       -\mathbb P(X^H_{g,n}=k)\right|
 \le 1-\frac{H^{\taut}_{g,n}(1)}{H_{g,n}(1)}
 =O(n^{-1}),
\end{equation}
and
\begin{equation}\label{eq:tautological-clt}
 \sup_{x\in\mathbb R}\left|
 \mathbb P\left(\frac{X^{\taut}_{g,n}-d_n/2}{\sqrt{\nu n}}\le x\right)
 -\Phi(x)\right|=O(n^{-1/2}).
\end{equation}
\end{theorem}

For the upper bound, we observe that the even cohomology classes vanishing on the preimage $U_n$ of the complement of the boundary strata of dimension at most one in $\Mbar_g$ are tautological.

\begin{lemma}\label{lem:tautological-support}
Fix $g\ge2$, and let $B_{\le1}\subset\Mbar_g$ be the union of the boundary strata $\mathcal M_\Gamma$ of dimension at most one. For the forgetful morphism $\pi_n:\Mbar_{g,n}\to\Mbar_g$, let $U_n=\pi_n^{-1}(\Mbar_g\setminus B_{\le1})$. Then, for every $k$,
\[
 \ker\bigl(H^{2k}(\Mbar_{g,n};\mathbb Q)\longrightarrow H^{2k}(U_n;\mathbb Q)\bigr)
 \subset RH^{2k}(\Mbar_{g,n}).
\]
Consequently, writing $b_c(Y)=\sum_j\dim H_c^j(Y;\mathbb Q)$,
\begin{equation}\label{eq:tautological-complement-bound}
 H_{g,n}(1)-H^{\taut}_{g,n}(1)\le b_c(U_n).
\end{equation}
\end{lemma}

\begin{proof}
The union $B_{\le1}$ is closed, since the boundary of a stratum consists of strata of smaller dimension. By \eqref{eq:stratum-dimension}, a stratum $\mathcal M_\Gamma$ in $B_{\le1}$ has all vertices of genus zero, except possibly one genus-one vertex of valency one. Forgetting markings contracts only rational components, so $Z_n=\pi_n^{-1}(B_{\le1})$ is covered by the images of finitely many gluing morphisms
\[
 \prod_v\Mbar_{g_v,n_v}\longrightarrow Z_n,
\]
one for each marked stable graph whose unmarked stabilization lies in $B_{\le1}$, where $n_v$ counts the flags and legs at $v$. Each source is a product of genus-zero factors and at most one genus-one factor.

Although $Z_n$ may be singular, these covers suffice. If $f:Y\to Z$ is a proper surjection with $Y$ smooth and proper, then
\begin{equation}\label{eq:pure-homology-image}
 \operatorname{im}\bigl(H_m(Y;\mathbb Q)\longrightarrow H_m(Z;\mathbb Q)\bigr)
 =W_{-m}H_m(Z;\mathbb Q).
\end{equation}
This is the dual form of Deligne's result \cite[Proposition~8.2.5]{DeligneHodgeIII}. Indeed, cohomological descent for the \v{C}ech nerve $Y_\bullet$ of $f$ gives a spectral sequence $E_1^{p,q}=H^q(Y_p;\mathbb Q)\Rightarrow H^{p+q}(Z;\mathbb Q)$ \cite[Section~8]{DeligneHodgeIII}. Each $Y_p$ is proper, so $H^q(Y_p)$ has weights at most $q$. In total degree $m$, the terms with $p\ge1$ therefore have weights at most $m-1$, whereas $H^m(Y)$ is pure of weight $m$. Hence the kernel of $H^m(Z)\to H^m(Y)$ is $W_{m-1}H^m(Z)$, and dualizing gives \eqref{eq:pure-homology-image}. The same argument applies to proper Deligne--Mumford stacks with rational coefficients.

Let $d=\dim\Mbar_{g,n}$. By localization and Poincar\'e duality, the kernel of restriction to $U_n$ in degree $2k$ is the image of
\[
 H_{2d-2k}(Z_n;\mathbb Q)(-d)\longrightarrow H^{2k}(\Mbar_{g,n};\mathbb Q).
\]
The target is pure of weight $2k$, so by strictness of the weight filtration and \eqref{eq:pure-homology-image}, applied to the disjoint union of the gluing sources, this image is spanned by Gysin pushforwards from the sources. These pushforwards shift degrees by even integers, so only the even cohomology of the sources is involved. That cohomology is tautological by Keel's theorem in genus zero \cite{Keel}, Petersen's theorem in genus one \cite{Petersen} and the K\"unneth formula, because odd cohomology occurs in at most one factor. Since gluing pushforwards preserve tautological classes, the inclusion follows. The even part of the quotient by the tautological classes therefore has dimension at most $\sum_kb_{2k}(U_n)\le b_c(U_n)$, by Poincar\'e duality on the smooth open stack $U_n$.
\end{proof}

To bound $b_c(U_n)$, we use the blowup construction in \cite{Kiem}.

\begin{lemma}\label{lem:relative-betti-substitution}
Fix $g\ge2$ and an open substack $V\subset\Mbar_g$. For $r,s\ge0$, let $\Mbar_{g,r|s}$ denote the Hassett space with $r$ markings of weight $1$ and $s$ markings of a positive rational weight $\epsilon$ satisfying $s\epsilon<1$. Let
\[
 U_{r|s}=\Mbar_{g,r|s}\times_{\Mbar_g}V,
 \qquad
 F_r(w)=\sum_{s\ge0}b_c(U_{r|s})\frac{w^s}{s!}.
\]
Then, as formal power series,
\begin{equation}\label{eq:relative-betti-substitution}
 \sum_{n\ge0}b_c(U_{n|0})\frac{z^n}{n!}
 =F_0(T(z,1)).
\end{equation}
\end{lemma}

\begin{proof}
By \cite[Section~3]{Kiem}, for $r\ge 1$, the natural morphism $\Mbar_{g,r|s}\to\Mbar_{g,r-1|s+1}$
is the composition of a sequence of blowups whose centers are disjoint unions of $\binom{s}{j}$ copies of $\Mbar_{g,r|s-j}$ for $2\le j\le s$.
By the blowup formula,
\[
 b_c(U_{r|s})
 =b_c(U_{r-1|s+1})
  +\sum_{j=2}^s\binom{s}{j}(j-1)\,b_c(U_{r|s-j}),\quad \text{for }r\ge 1
\]
that is,
\[
 F_r(w)=F'_{r-1}(w)
       +F_r(w)\sum_{j\ge2}(j-1)\frac{w^j}{j!},
 \qquad\text{or}\qquad
 F_r(w)=\frac{F'_{r-1}(w)}{R_w(w,1)},
\]
since $R_w(w,1)=1-\sum_{j\ge2}(j-1)w^j/j!$ by \eqref{eq:inverse-definition}. Because $R(T(z,1),1)=z$, the chain rule gives $F_r(T(z,1))=\partial_zF_{r-1}(T(z,1))$, hence $F_r(T(z,1))=\partial_z^rF_0(T(z,1))$. Since $T(0,1)=0$ and $F_r(0)=b_c(U_{r|0})$, evaluating the preceding identity at $z=0$ gives
\[
 \left.\partial_z^r F_0(T(z,1))\right|_{z=0}
 =b_c(U_{r|0})
 \qquad (r\ge0).
\]
Thus the formal Taylor expansion of $F_0(T(z,1))$ at $z=0$ is
\[
 F_0(T(z,1))
 =\sum_{r\ge0}b_c(U_{r|0})\frac{z^r}{r!},
\]
which proves \eqref{eq:relative-betti-substitution}.
\end{proof}

\begin{proof}[Proof of Theorem~\ref{thm:tautological-dominance}] 

For the upper bound, we apply Lemma~\ref{lem:relative-betti-substitution} with $V=\Mbar_g\setminus B_{\le1}$. We first bound the light series $F_0(w)$. Every core stratum $\mathcal M_\Gamma$ contained in $V$ has dimension at least two, so \eqref{eq:stratum-dimension} gives
\[
3g-3-|E(\Gamma)|=\dim\mathcal M_\Gamma\ge2.
\]
Thus every such graph has at most $3g-5$ edges.

Fix one of these graphs $\Gamma$. As in the proof of Proposition~\ref{prop:graph-sum}, we stratify the space of light curves with core graph $\Gamma$ by assigning labels to the smooth locus of the stabe core and to the chains inserted at its nodes. For an assignment with $r$ labels on the smooth locus of the stable core and $k_e$ labels at each edge $e$, the corresponding stratum before taking the finite quotient is
\[
X_{\Gamma,r}\times\prod_{e\in E(\Gamma)}\mathrm{LM}_{k_e},\qquad r+\sum_e k_e=n.
\]
The exponential Betti bound \eqref{eq:exponential-betti-bound} gives $b_c(X_{\Gamma,r})\le K_\Gamma B_\Gamma^r$, with constants independent of $r$, while Lemma~\ref{lem:node-factors} at $t=1$ gives $b_c(\mathrm{LM}_k)=k!$, including $k=0$. Hence
\[
\sum_{r\ge0}b_c(X_{\Gamma,r})\frac{w^r}{r!}\preceq K_\Gamma e^{B_\Gamma w},\qquad \sum_{k\ge0}b_c(\mathrm{LM}_k)\frac{w^k}{k!}=\frac1{1-w},
\]
where $\preceq$ denotes coefficientwise inequality between power series with nonnegative coefficients.

There are $n!/(r!\prod_e k_e!)$ assignments with these values of $r$ and $k_e$. The K\"unneth formula and summation over all assignments therefore give the product of the preceding generating series, with one factor $1/(1-w)$ for each edge. Passing to a finite quotient can only decrease the compactly supported rational Betti numbers, and the total compactly supported Betti number of a space with a finite stratification is at most the sum over its strata. Summing over the finitely many isomorphism classes of core graphs occurring in $V$, we obtain constants $K_g,B_g>0$ such that
\begin{equation}\label{eq:light-complement-majorant}
F_0(w)\preceq\sum_\Gamma\frac{K_\Gamma e^{B_\Gamma w}}{(1-w)^{|E(\Gamma)|}}\preceq\frac{K_g e^{B_gw}}{(1-w)^{3g-5}}.
\end{equation}
The coefficients of $T(z,1)$ are nonnegative (Lemma~\ref{lem:inverse}), so Lemmas~\ref{lem:tautological-support} and~\ref{lem:relative-betti-substitution} give
\[
H_{g,n}(1)-H^{\taut}_{g,n}(1)\le K_g\,n!\,[z^n]\frac{e^{B_gT(z,1)}}{(1-T(z,1))^{3g-5}}.
\]
By \eqref{eq:T-expansion} at $t=1$, as $z\to e-2$ in the $\Delta$-domain of Lemma~\ref{lem:inverse},
\begin{equation}\label{eq:T-expansion-at-one}
1-T(z,1)=\sqrt{\frac{2(e-2)}{e}}\left(1-\frac{z}{e-2}\right)^{1/2}+O\!\left(1-\frac{z}{e-2}\right).
\end{equation}
Consequently,
\[
\frac{e^{B_gT(z,1)}}{(1-T(z,1))^{3g-5}}=e^{B_g}\left(\frac{e}{2(e-2)}\right)^{(3g-5)/2}\left(1-\frac{z}{e-2}\right)^{-(3g-5)/2}\left(1+O\!\left(\left(1-\frac{z}{e-2}\right)^{1/2}\right)\right).
\]
Applying the transfer theorem gives
\[
[z^n]\frac{e^{B_gT(z,1)}}{(1-T(z,1))^{3g-5}}=O\!\left((e-2)^{-n}n^{(3g-5)/2-1}\right).
\]
Combining this estimate with the preceding coefficientwise bound yields
\begin{equation}\label{eq:non-tautological-total-bound}
H_{g,n}(1)-H^{\taut}_{g,n}(1)=O\!\left(n!\,(e-2)^{-n}n^{(3g-7)/2}\right).
\end{equation}
On the other hand, \eqref{eq:coefficient-asymptotics} at $t=1$ gives
\[
H_{g,n}(1)=C_g(1)\,n!\,(e-2)^{-n}n^{(3g-5)/2}\bigl(1+O(n^{-1/2})\bigr),\qquad C_g(1)>0.
\]
Dividing the two estimates, we obtain
\[
1-\frac{H^{\taut}_{g,n}(1)}{H_{g,n}(1)}=\frac{H_{g,n}(1)-H^{\taut}_{g,n}(1)}{H_{g,n}(1)}=O(n^{-1}),
\]
which proves the upper bound in \eqref{eq:tautological-fraction}.

For the lower bound, substitute $w=T(z,1)$ in \eqref{eq:fixed-hodge-difference-pole}. Using \eqref{eq:T-expansion-at-one} and applying the transfer theorem as above, we obtain
\begin{equation}\label{eq:fixed-hodge-difference-asymptotic}
\sum_{p-q=22}h^{p,q}(\Mbar_{g,n})=\frac{a_g}{\Gamma((3g-5)/2)}\left(\frac{e}{2(e-2)}\right)^{(3g-5)/2}n!\,(e-2)^{-n}n^{(3g-7)/2}\bigl(1+O(n^{-1/2})\bigr).
\end{equation}
Every class counted in \eqref{eq:fixed-hodge-difference-asymptotic} has even degree and lies off the Hodge diagonal, whereas tautological classes are algebraic and hence of type $(k,k)$. Therefore
\[
H_{g,n}(1)-H^{\taut}_{g,n}(1)\ge\sum_{p-q=22}h^{p,q}(\Mbar_{g,n}).
\]
Since $a_g>0$, dividing by \eqref{eq:coefficient-asymptotics} at $t=1$ proves the lower bound in \eqref{eq:tautological-fraction}.

The total variation bound \eqref{eq:tautological-tv} follows immediately from the degreewise inclusion of tautological cohomology and \eqref{eq:tautological-fraction}. Combining this bound with the CLT for $X^H_{g,n}$ in Theorem~\ref{thm:clt} gives \eqref{eq:tautological-clt}.
\end{proof}

By \eqref{eq:full-even-ratios}, the tautological classes also account for a proportion $1-O(n^{-1/2})$ of the full cohomology of $\Mbar_{g,n}$. The classes used in the lower bound are not algebraic. Non-tautological algebraic classes also exist, as the following example shows.

\begin{remark}\label{rem:genus-two-n20}
By \cite[Remark~3.10]{PetersenGenusTwo}, the quotient of the even cohomology of $\Mbar_{2,20}$ by its tautological part is concentrated in degree $22$ and has Hodge type $(11.11)$. As an $S_{20}$-representation it is the irreducible representation associated with the partition $(2^{10})$, which is of dimension $20!/(10!\,11!)=16796$ by the hook-length formula. Hence
\[
 H^{\taut}_{2,20}(t)=H_{2,20}(t)-16796\,t^{11}.
\]
For instance, the non-tautological algebraic class on $\Mbar_{2,20}$ constructed by Graber and Pandharipande \cite{GraberPandharipande} lies in degree $22$ and has nonzero image in this quotient.
\end{remark}

\providecommand{\bysame}{\leavevmode\hbox to3em{\hrulefill}\thinspace}

\end{document}